\documentclass[11pt]{amsart}

\usepackage{graphicx}
\usepackage[dvips]{color}
\usepackage[top=26mm, bottom=26mm, left=32mm, right=32mm]{geometry}

\usepackage{amssymb,amsmath,bm}
\usepackage{xcolor}
\usepackage{lineno}
\usepackage[colorlinks]{hyperref}
\usepackage{graphicx} 
\usepackage{subcaption}
\usepackage{csvsimple}
\usepackage{multirow}
\usepackage{array}
\usepackage{algorithm}
\usepackage{algpseudocode}

\renewcommand{\d}{{\rm{d}}}
\newcommand{\vct}[1]{\bm{#1}}
\newcommand{\mtx}[1]{\mathbf{#1}}

\newcommand{\vtwo}[2]{\left[\begin{array}{c} #1 \\ #2 \end{array}\right]}
\newcommand{\vthree}[3]{\left[\begin{array}{c} #1 \\ #2 \\ #3 \end{array}\right]}

\newcommand{\mtwo}[4]{\left[\begin{array}{cc}    #1 & #2 \\ #3 & #4  \end{array}\right]}
\newcommand{\mthree}[9]{\left[\begin{array}{ccc} #1 & #2 & #3 \\
                                                 #4 & #5 & #6 \\
                                                 #7 & #8 & #9 \end{array}\right]}

\newtheorem{remark}{Remark}
\numberwithin{remark}{section}
\newtheorem{definition}{Definition}
\numberwithin{definition}{section}

\begin{document}
\title{An efficient boundary integral equation solution technique for solving aperiodic scattering problems from two-dimensional, periodic boundaries}
\author{R.Fisher$^1$, F. Agocs$^2$, and A. Gillman$^1$}
 \address{$1$ Department of Applied Mathematics, University of Colorado, Boulder}
 \address{$2$ Department of Computer Science, University of Colorado, Boulder}

\begin{abstract}
This manuscript presents an efficient boundary integral equation technique for solving two-dimensional Helmholtz problems defined in the half-plane bounded by an infinite, periodic curve with Neumann boundary conditions and with an aperiodic point source. 
The technique is designed for boundaries where one period does not require a large number of discretization points to achieve high accuracy. The inverse Floquet--Bloch transform turns the problem into evaluating a contour integral where the integrand is the solution of quasiperiodic boundary value problems. To approximate the integral, one must solve a collection of these problems. This manuscript uses a variant of the periodizing scheme by Cho and Barnett which avoids using the quasiperiodic Green's function in the integral representation and is amenable to a large amount of precomputation that can be reused for all the necessary solves. The solution technique is accelerated using low rank linear algebra and corner compression. The performance of the method is measured on a range of simple and complex boundary geometries, including in the presence of one or more corners with straight or curved edges present in the unit cell. The numerical results illustrate that the presented method is competitive with techniques using the quasiperiodic Green's function kernel for smooth boundaries, and up to 20--30 faster when corner compression and high accuracy is required.
\end{abstract}

\maketitle

\section{Introduction}
\label{sec:intro}

Periodic structures are commonly used in engineering applications to guide or localize electromagnetic and acoustic waves. They were first rigorously studied for the purposes of manufacturing and optimizing diffraction gratings \cite{mciver1998construction}, where they have been observed to give rise to modes that travel along the surface without attenuation. 
These guided (``trapped") surface waves
are defined as localized solutions of the wave equation without sources. They exist at a well-defined frequency, have finite energy, and are characterized by exponential decay in at least one direction. In the context of acoustic waves near periodic, open geometries, such as gratings, they are referred to as Rayleigh--Bloch surface waves, and have been studied extensively (see, e.g.\ \cite{PorterEvans1999Rayleigh,pagneux2013trapped} and many references therein). 
Trapped modes have been exploited for filtering and sensing, such as decreasing the impact of water waves on offshore structures \cite{mciver2002offshore,Kagemoto1992inverse,cobelli2011experimental}, sound insulation or manipulation \cite{guancongacousticmetareview,Li2014acousticmeta}, and in phononic and photonic crystals \cite{Khelif2003ultrasonic}. 

The existence of unattenuated surface waves, however, poses a challenge for numerical simulation, since truncation of the otherwise unbounded computational domain results in large reflections and numerical errors near the boundary. 
If the wave incident on the periodic boundary is quasiperiodic, i.e.\ periodic up to a complex phase, such as a plane wave, the solution is also expected to be quasiperiodic and may therefore be periodized (computed over one unit cell of the periodic domain). 
In the case of an aperiodic incident field (e.g.\ the subject of this work, a cylindrical wave emanating from a line source), however, the solution does not possess any symmetry that may be exploited to make the computation tractable.
Such aperiodic incident fields may still be considered via the Floquet--Bloch transform, a generalization of the Fourier transform, which decomposes them into quasiperiodic components. Each of the resulting fields defines a quasiperiodic scattering problem, parameterized by a Bloch phase. These individual problems may be solved by periodization. 
An integral of the quasiperiodic solutions over the Bloch phase then needs to be evaluated.

The Floquet--Bloch transform has been used to solve aperiodic scattering problems from 2D periodic arrays of Dirichlet obstacles \cite{wu1966properties,capolino2005mode}, continuous Neumann boundaries with corners \cite{agocs2024trapped}, periodic closed waveguides \cite{zhang2021numerical}, and locally perturbed periodic surfaces \cite{zhang2018perturbed}. Three-dimensional, bi-periodic surfaces with Dirichlet or impedance boundary conditions have also been investigated \cite{lechleiter2017non,arens2024high}, where the Floquet--Bloch transform involves a nested integral. 

Assuming time-harmonicity, the acoustic scattering solution satisfies the two-dimensional Helmholtz equation in the upper half plane where the bottom boundary is defined by an infinite periodic curve $\Gamma$.
Let $\Omega$ denote the domain where the partial differential equation is defined, as pictured in Figure~\ref{fig:unitcell}(a). Throughout the manuscript, $d$ denotes the length of one period of $\Gamma$.
The Helmholtz problem for incident frequency $\omega$ and unknown total field $u_{\mathrm{tot}}(\vct{x})$ is then
\begin{equation}
    - (\Delta + \omega^2)u_{\mathrm{tot}} = \delta(\vct{x} - \vct{x}_0) \quad \text{in }\Omega,
\end{equation}
where $\vct{x}_0$ is the location of a point source. The boundary is a sound-hard surface resulting in Neumann boundary conditions
\begin{equation}
    (u_{\mathrm{tot}})_{\vct{\nu}} = \vct{\nu}\cdot \nabla u_{\mathrm{tot}} = 0 \quad \text{on }\Gamma,
\end{equation}
where $ \vct{\nu}$ is the upwards-pointing unit normal vector at $\vct{x}\in \Gamma$. The total field is defined as 
$u_{\rm tot} = u_i + u,$ where $u_i$ is a known incident wave generated by the point source at $\vct{x}_0$ and $u$ is the scattered field defined as the solution to the boundary value problem
\begin{equation}\label{eq:bvp}
\begin{aligned}
- (\Delta + \omega^2)u &= 0 \quad &&\text{in }\Omega,\\    
    u_{\vct{\nu}} &= -(u_i)_{\vct{\nu}} \quad &&\text{on }\Gamma,
    \end{aligned}
\end{equation}
with a radiation condition in the $y$-direction. 

By applying the inverse Floquet--Bloch transform to \eqref{eq:bvp}, we find that the solution can be 
found by carrying out the Floquet--Bloch integral
\begin{equation}
    u = \frac{d}{2\pi} \int_{-\frac{\pi}{d}}^{\frac{\pi}{d}} u_{\kappa}\d \kappa,
    \label{eq:fb-integral}
\end{equation}
where 
$u_\kappa$ is the solution of 
\begin{equation}\label{eq:quasi-prob}
\begin{aligned}
    -(\Delta + \omega^2)u_{\kappa} &= 0 \quad &&\text{in } \Omega, \\
    (u_\kappa)_{\vct{\nu}} &= -f_{\vct{\nu}} \quad &&\text{on } \Gamma. 
\end{aligned}
\end{equation}
$f_{\vct{\nu}}$ is the normal derivative of the incident field generated by a quasiperiodic array of point sources, referred to as the quasiperiodic Green's function,
\begin{equation}\label{eq:qp-rhs}
    f( \vct{x}) = \frac{i}{4}\sum_{l \in \mathbb{Z}} \alpha^l H_0^{(1)}(| \vct{x} -  \vct{x}_0 - l \vct{d}|), 
\end{equation}
with $ \vct{d} = (d, 0)$. The parameter $\alpha = e^{i\kappa d}$ denotes the Bloch phase. Thus, for each $\kappa$ in (\ref{eq:fb-integral}) there is a different Bloch phase $\alpha$.
 The radiation condition in the $y$-direction can be defined via a Rayleigh--Bloch expansion. This is a Fourier series with discrete horizontal and vertical wavenumbers, defined as 
\begin{equation}\label{eq:ray}
    u_\kappa (x,y) = \sum_{m \in \mathbb{Z}} b_me^{i(\beta_m x + k_my)},
\end{equation}
where
\begin{equation}\label{eq:wavenumbers}%
\begin{aligned}
    \beta_m &= \kappa + 2\pi m/d, \\
    k_m &= 
    \begin{cases}
        \sqrt{\omega^2 - \beta_m^2} \quad &\text{if } \omega \geq |\beta_m|, \\ 
        i \sqrt{\beta_m^2 - \omega^2} \quad &\text{otherwise,}
    \end{cases}
\end{aligned}
\end{equation}
for $m \in \mathbb{Z}$, and $b_m$ are unknown expansion coefficients. The solution $u_\kappa$ to this half-space boundary value problem satisfies the following quasiperiodicity
condition:

\begin{equation}
\label{eq:qpcond}
u_\kappa( \vct{x}) = \alpha^{-1} u_\kappa( \vct{x}+\vct{d}).
\end{equation}
Algorithm~\ref{alg:high-level} provides an overview of the proposed solution technique for evaluating the integral \eqref{eq:fb-integral}.

There are a variety of methods for approximating solutions to the quasiperiodic scattering problem (\ref{eq:quasi-prob}). For the purposes of this manuscript, we focus on methods designed for recasting (\ref{eq:bvp}) as an integral equation. 
A commonly used approach is to represent the solution using the quasiperiodic Green's function \cite{desanto1998theoretical,bruno2009efficient,pinto2021fast,meng2023new,agocs2024trapped}. One of the main challenges with this approach is that the quasiperiodic Green's function is defined as an infinite sum that can be slowly convergent, and is not convergent at Wood's anomalies. There exist fast algorithm techniques for handling the infinite sum, many of which are reviewed in \cite{moroz2006quasi}. Another drawback of evaluating the quasiperiodic Green's function directly is that it needs to be recomputed for each Bloch phase $\alpha$. When using the quasiperiodic Green's function in the integral representation of the solution, a new quasiperiodic Green's function must be computed for each discretization point on the boundary. This means that in the approximation of (\ref{eq:fb-integral}), the quasiperiodic Green's function needs to be evaluated many times. One example is \cite{agocs2024trapped}, which uses the lattice sum approach from \cite{yasumoto2002efficient} to compute the quasiperiodic Green's function. It then defines and solves a new integral equation at each Bloch phase involved in the inverse Floquet--Bloch transform, recomputing all matrices each time. Due to the breakdown of the lattice sum approach, the Floquet--Bloch integral involved needs to avoid Bloch phases where Wood's anomalies occur.

An alternative approach for quasiperiodic scattering problems that avoids the use of the quasiperiodic Green's function in the layer representations was first presented in \cite{BARNETTgreengard}. It was further generalized in \cite{Cho:15}. The idea is to represent the solution via integral operators that involve the free-space Green's function and a linear combination of basis functions with unknown coefficients. The unknowns are found by solving a system found by enforcing the boundary condition and quasiperiodicity. The quasiperiodic scattering problems considered in this manuscript still require the quasiperiodic Green's function for evaluating the boundary data. The quasiperiodic Green's function is represented in a similar manner as the integral representation, which avoids the infinite sum. The main advantages to this approach are that 1) it is robust at Wood's anomalies, and 2) the system is separable into Bloch phase-dependent and independent pieces. This allows a significant portion of a direct solver to be precomputed and used for all solves independent of the Bloch phase. In \cite{2013_Gillman}, exploiting this precomputation led to a 600 times speedup over solving the same linear system using an FMM-accelerated GMRES solver for a problem involving over 200 right-hand-sides. The ideas in \cite{2013_Gillman} have been extended to multilayered media scattering problems involving periodic interfaces \cite{Zhang21,zhang2022,ChoWu}. These fast solvers are only necessary if the number of discretization points needed to resolve one period of the boundary is large.

For the problems of interest in this manuscript, not many points are needed to discretize the boundary and obtain accurate approximate solutions to equation (\ref{eq:bvp}). Thus, we choose to use the periodizing method presented in \cite{Cho:15} modified for a Neumann boundary condition on the interface. We created a new efficient direct solver that is inspired by the methods in \cite{2013_Gillman,Zhang21,zhang2022} for the linear system that results from the discretization of the integral equation. This solver is specifically designed for problems where the system is not large enough to justify the use of the fast direct solvers from \cite{2013_Gillman,Zhang21,zhang2022}. For the systems under consideration, the solver has a relatively lean precomputation, reuses information that is independent of the Bloch phase, and can be applied  quickly for the multiple solves needed in the Floquet--Bloch integral.

The goal of this manuscript is to illustrate how several numerical methods---including the Floquet--Bloch approach, a boundary integral equation method, an associated fast solver, and the recycling of matrices involved---may be combined to significantly reduce the computation time needed to solve aperiodic scattering problems in periodic geometries, by pushing more of the computational effort into operations carried out once at the start of the process. 

There are several novel contributions presented in this manuscript. It is the first presentation of localized (or corner) compression for a problem not involving a single geometry in free space. This is also the first use of directional compression techniques in the quasiperiodic setting. The introduction of this technique reduces the size of the linear system that needs to be solved and the cost of constructing low rank factorizations. More generally, this work illustrates the advantages of using compression techniques in situations where the system size is not large enough to justify the use of a \emph{fast} direct solver. The efficiency of the solver allows for exploration of the dispersion relation of trapped modes of a novel surface (see Figure \ref{fig:fb-integral}), and scattering from three-dimensional surfaces.

The manuscript proceeds by first introducing the Floquet--Bloch integral in Section \ref{sec:flobloch}. The boundary integral periodization technique is presented in Section \ref{sec:boundInt}. This section also includes a technique for representing the two-dimensional, half-space quasiperiodic Helmholtz Green's function via the free-space Helmholtz Green's function. Section \ref{sec:solver} presents the proposed efficient direct solver technique including a method of handling boundaries with corners. Section \ref{sec:numerics} illustrates the performance of the proposed solver for solving the boundary value problem in (\ref{eq:bvp}) and computing the Floquet--Bloch integral. Section \ref{sec:conclusion} closes the document with a summary and a discussion of the goal to extend this technique to three-dimensional problems.

\begin{algorithm}[ht]
\caption{Overview of the aperiodic solution technique}
\label{alg:high-level}
	\renewcommand{\algorithmicrequire}{\textbf{Input:}}
	\renewcommand{\algorithmicensure}{\textbf{Output:}}
\begin{algorithmic}[1]
\Require aperiodic point source location $\vct{x}_0$ and 
    two-dimensional periodic boundary $\Gamma$.
        \Ensure  approximation of $u = \dfrac{d}{2\pi}\displaystyle\int_{-\pi/d}^{\pi/d} u_\kappa\,\d\kappa$.

    \State Create $N_\kappa$ quadrature nodes $\{\kappa_j\}$ and weights $\{w_j\}$ on the integration contour $\kappa(s)$ defined in (\ref{eq:contour}).

    \State Precompute the $\alpha$-independent portions of the Helmholtz direct solver as described in Algorithm \ref{alg:precomp}.

    \For{$j = 1,\ldots,N_\kappa$}
        \State Create the boundary data for $\kappa_j$ via  the method in Section \ref{sec:gqp}.
        
        \State Solve a quasiperiodic scattering problem to find $u_{\kappa_j}$ using Algorithm \ref{alg:solve}.
    
    \EndFor

    \State Evaluate
    $
        u \gets \dfrac{d}{2\pi}
        \displaystyle\sum_{j=1}^{N_\kappa} w_j u_{\kappa_j}
    $.
    \State \Return $u$
\end{algorithmic}
\end{algorithm}

\section{Floquet--Bloch integral }
\label{sec:flobloch}

A common way of making the aperiodic problem tractable is to represent the scattered solution $u$ as an integral of quasiperiodic solutions via the inverse Floquet--Bloch transform. The transform utilizes the identity
\begin{equation}
\begin{aligned}
\label{eq:source}
    \delta(\vct{x} - \vct{x}_0) &= \frac{d}{2\pi} \int_{-\frac{\pi}{d}}^{\frac{\pi}{d}} \sum_{n \in \mathbb{Z}} \delta(\vct{x} -\vct{x}_0 - n\vct{d})e^{i n \kappa d}\d \kappa\\
    &= \frac{d}{2\pi} \int_{-\frac{\pi}{d}}^{\frac{\pi}{d}} \sum_{n \in \mathbb{Z}} \delta(\vct{x} -\vct{x}_0 - n\vct{d})\alpha^n\d \kappa.
\end{aligned}
\end{equation}
The left-hand-side of \eqref{eq:source} is the common representation of a point source (the delta function). The integrand on the right-hand-side of \eqref{eq:source} is an infinite, quasiperiodic array of point sources with an associated constant phase difference $\alpha$ from Section \ref{sec:intro}.
By linearity, it follows that the aperiodic scattered solution $u$ to equation~\eqref{eq:bvp} may be written as the integral in equation \eqref{eq:fb-integral}, and that $u_{\kappa}$ satisfies equation (\ref{eq:quasi-prob}).
The quasiperiodic scattered solution $u_\kappa$ is unique at all values of $\kappa \in [-\pi/d, \pi/d)$, except potentially a discrete set \cite{bonnet1994guided}. These values of $\kappa$ correspond to on-surface wavenumbers of trapped modes, which do not decay in the horizontal direction but decay exponentially in the vertical direction.

A further consequence of quasiperiodicity is that on-surface wavenumbers are only unique modulo $2\pi/d$, since $\kappa$ and $\beta_n = \kappa + 2n\pi/d$ for $n \in \mathbb{Z}$ give the same Bloch phase $\alpha$. For an incident frequency $\omega$ and real on-surface wavenumber $\beta_n$, the vertical wavenumber $k_n$ is then either positive real or positive imaginary, as given by \eqref{eq:wavenumbers}.
Physically, having $k_n$ be purely real means the vertical component of the solution is a propagating (outgoing) wave, whereas it decays exponentially for a purely imaginary $k_n$. The $k_n = 0$ case, where the behavior of the vertical component of the scattered field abruptly changes, corresponds to a Wood's anomaly \cite{wood1902}.

\begin{figure}
\centering
\begin{subfigure}[t]{0.445\textwidth}
    \centering
    \includegraphics{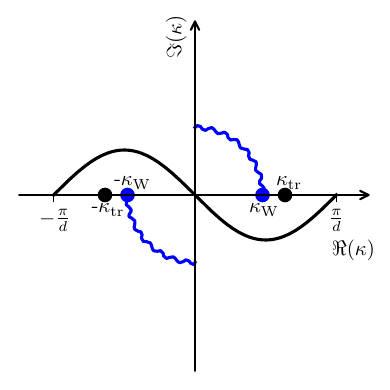}
    \caption{Sketch of the Floquet--Bloch integration contour with branch cuts and poles labeled.}
\end{subfigure}
\hfill
\begin{subfigure}[t]{0.445\textwidth}
    \centering
    \includegraphics{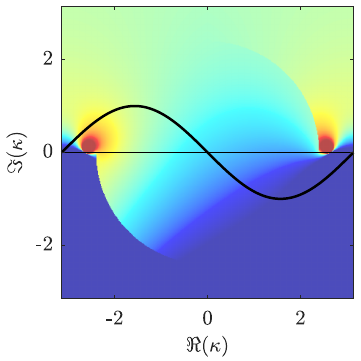}
    \caption{Real part of the integrand of the inverse Floquet--Bloch transform for the staircase geometry.}
\end{subfigure}
\begin{subfigure}[t]{0.085\textwidth}
    \centering
    \includegraphics{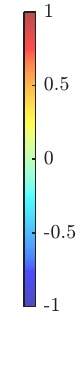}
\end{subfigure}
\begin{subfigure}[t]{0.445\textwidth}
    \centering
    \includegraphics{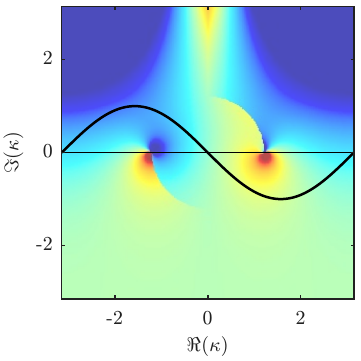}
    \caption{Real part of the integrand of the inverse Floquet--Bloch transform for the dinosaur geometry.}
\end{subfigure}
\hfill
\begin{subfigure}[t]{0.445\textwidth}
    \centering
    \includegraphics{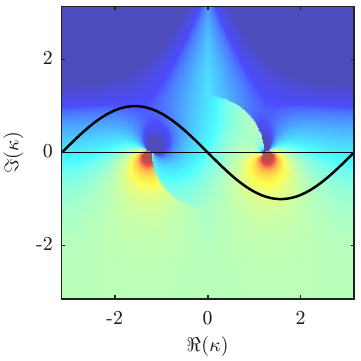}
    \caption{Real part of the integrand of the inverse Floquet--Bloch transform for the wave geometry.}
\end{subfigure}
\hfill
\begin{subfigure}[t]{0.085\textwidth}
    \centering
    \includegraphics{fb_colorbar.pdf}
\end{subfigure}
\caption{\label{fig:fb-integral} Deformation of the integral contour for the inverse Floquet--Bloch transform and real part of its integrand. (a) The black, sinusoidal curve shows the deformed contour that goes around the poles corresponding to trapped modes ($\kappa_{\mathrm{tr}}$, black dots) and branch cuts corresponding to Wood's anomalies ($\kappa_{\mathrm{W}}$, blue dots and wavy lines). (b-d) Show the real part of the integrand $u_{\kappa}$ in the inverse Floquet--Bloch transform \eqref{eq:fb-integral}, with the branch cuts and poles visible and the contour deformation in solid black, for various geometries defined in section \ref{sec:numerics}. The source and target points used for generating these figures are as follows. (b) Staircase geometry from the source point $\vct{x}_0 = (-0.2, 0.6)$, evaluated at the target point $\vct{x} = (0.22, 0.34)$, with $\omega = 2.4$. (c) Dinosaur geometry from the source point $\vct{x}_0 = (0.4, 0.45)$, evaluated at the target point $\vct{x} = (-0.35, 0.4)$, with $\omega = 1.2$. (d) Wave geometry from the source point $\vct{x}_0 = (0.4, 0.45)$, evaluated at the target point $\vct{x} = (-0.35, 0.4)$, with $\omega = 1.2$.}
\end{figure}

Figure~\ref{fig:fb-integral}(a) is a sketch of the integrand $u_{\kappa}$ of the Floquet--Bloch transform in the complex $\kappa$-plane. The values $\kappa_{\mathrm{tr}}$ denote the $\kappa$ values that correspond to trapped modes, which present as poles. The Wood's anomalies, denoted by $\kappa_W$, can be shown to be $z^{-1/2}$-type singularities \cite{BARNETTgreengard}. Both $\kappa_{\mathrm{tr}}$ and $\kappa_W$ are located on the real axis. There are branch cuts emanating from the Wood's anomalies, the directions and shapes of which are determined by how $k_n$ is evaluated for a general complex-valued $\kappa$ (and $\beta_n$). In this work, the choice in \eqref{eq:wavenumbers} leads to branch cuts that are circular arcs, and together with the choice of integration contour for the Floquet--Bloch integral, ensures that the solution is outgoing.

Figure \ref{fig:fb-integral}(b--d) illustrate the real part of the integrand $u_\kappa$ for various geometries (defined in section \ref{sec:numerics}) and wavenumbers. The poles on the real axis at $\kappa = \kappa_{\mathrm{tr}}$ are visible, as are the Wood's anomalies and their arc-shaped branch cuts. While the integrand changes for the different geometries, the locations of the poles and Wood's anomalies as well as the shape and location of the branch cuts remain consistent for a fixed wavenumber (see Figure \ref{fig:fb-integral}(c-d)). 

To evaluate the integral in equation (\ref{eq:fb-integral}), we define a contour parameterized by $s = \Re(\kappa)$ that avoids the branch cuts and the Wood's anomalies. In Figure~\ref{fig:fb-integral} and the numerical examples in this manuscript, the parameterization is given by
\begin{equation} \label{eq:contour}
    \kappa(s) = s - i\sin(sd), \quad s \in \bigg[ -\frac{\pi}{d}, \frac{\pi}{d}\bigg). 
\end{equation}
This choice of contour always avoids the singularities present on the real axis and is therefore fixed during our numerical experiments and for all values of $\omega$ (see Figure \ref{fig:fb-integral}). In accordance with the limiting absorption principle \cite{sveshnikov1950radiation}, the contour is required to pass through the second and fourth quadrants of the complex plane to ensure the numerical solution is outgoing. Existing choices of the deformed contour involve piecewise smooth \cite{zhang2021numerical,bakhtafrouz2015application}, $\mathbb{C}^2$ \cite{shanin2024double}, and analytic \cite{agocs2024trapped} functions. We choose the latter option of an analytic contour, whose amplitude (extent of deformation) was set by checking how accurately a Floquet--Bloch integral along it can reconstruct the free-space, two-dimensional Helmholtz Green's function.

To approximate the integral in equation (\ref{eq:fb-integral}), quadrature points are placed on $\kappa(s)$. Since the integrand is periodic (the points $\kappa = \pm \pi/d$ are equivalent), the periodic trapezoidal rule is quickly convergent \cite{ptrconvergence}. However, in the limits $\omega \to 0$ and $\omega \to \omega_c$, where $\omega_c$ is the largest frequency at which trapped modes exist, the integration contour $\kappa(s)$ passes between two poles coalescing at $\kappa_{\mathrm{tr}} = 0$ and $\kappa_{\mathrm{tr}} = \pm \pi/d$, respectively. For example, Figure \ref{fig:fb-integral}(b) shows an example of this for the staircase geometry with $\omega=2.4$. Here, the poles at $\kappa_{\mathrm{tr}}$ approach $\pm \pi/d$. In these cases, we increase the density of quadrature nodes in the vicinity of the coalescing poles according to \cite{barnett2018unified,agocs2024trapped}. The method works by taking quadrature nodes that are evenly spaced in $s$ along the path \eqref{eq:contour}, and modifying their density to $\theta(s)$, a function numerically computed from its derivative that takes its maximum near the point(s) where bunching is required. The function $\theta(s)$ is defined as
\begin{equation}\label{eq:expgrading}
    \kappa(s) = \theta(s) - i\sin(\theta(s)), \quad \theta'(s) = \begin{cases}
        \alpha \cosh(b\sin \tfrac{sd}{2}) \; &\text{when }\kappa_{\mathrm{tr}} \to \pm \pi/d, \\ 
        \alpha \cosh(b\sin \tfrac{sd - \pi}{2}) \; &\text{when }\kappa_{\mathrm{tr}} \to 0 \\
    \end{cases}
\end{equation}
where $b$ is a free parameter, the normalization parameter $\alpha$ is obtained from the requirement that $\int_0^{2\pi/d} \theta'(s) \d s = 2\pi/d$, and $\theta(s)$ is computed from $\theta'(s)$ by taking the antiderivative of its spectral interpolant using the FFT. Using \eqref{eq:expgrading} results in nodes that are evenly spaced in $s$ being mapped to be a factor $e^b$ closer together near $\kappa_{\mathrm{tr}} \in \{0, \pm \pi/d \}$. The value of $b$ can be adjusted as needed. 

At each quadrature node, the solution of one quasiperiodic scattering problem with the same wavenumber but different Bloch phase $\alpha$ and boundary condition is required. Thus it is necessary to have an efficient numerical method for approximating the value $u_\kappa$.

\section{The periodizing scheme}
\label{sec:boundInt}


The periodization scheme used in this manuscript is a variation of the method presented in \cite{Cho:15} that has been modified for an interface with Neumann boundary conditions. The method in \cite{Cho:15} can be viewed as a modified version of the integral formulation in \cite{1978_kress,rokh83}. The presentation of the periodization scheme is kept  relatively brief, since detailed explanations are available in \cite{Cho:15,Zhang21,zhang2022}. 

This section breaks the periodization scheme into three pieces: sections~\ref{sec:geom}--\ref{sec:lin-sys} describe the partitioning of the domain, the representation of the solution to the Neumann problem, and the linear system that must be solved, respectively. Section \ref{sec:gqp} presents the application of the periodization technique for representing the quasiperiodic Green's function. 

\subsection{Partitioning of the geometry}
\label{sec:geom}
The periodization scheme begins by partitioning the domain $\Omega$ and introducing the idea of a \emph{unit cell}. The unit cell is a rectangle where the left and right boundaries are defined to bound one copy of the periodic surface. The top boundary is placed above at a chosen location above the surface. The bottom boundary can lie anywhere below the periodic surface. Figure~\ref{fig:unitcell}(b) illustrates the unit cell. Let $x = x_L$ and $x = x_R$ denote the left and right bounds, and $y= y_U$ denote the upper bound of the unit cell. For the Neumann boundary value problem, the periodic surface is the lower boundary of the unit cell. Let $\Gamma_0$ denote the portion of the lower boundary inside the unit cell. Similarly, $\Gamma_{-1}$ and $\Gamma_1$ denote the portions of the lower boundary in the left and right neighboring copies of $\Gamma_0$, respectively. 

\subsection{Representation of the solution}
\label{sec:rep-soln}
\begin{figure}
\centering
\begin{subfigure}[t]{0.49\textwidth}
    \centering
    \includegraphics{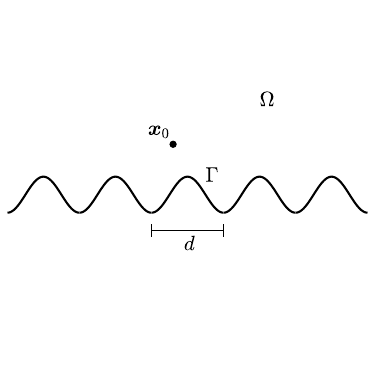}
    \caption{Domain and boundary.}
\end{subfigure}
\hfill
\begin{subfigure}[t]{0.49\textwidth}
    \centering
    \includegraphics{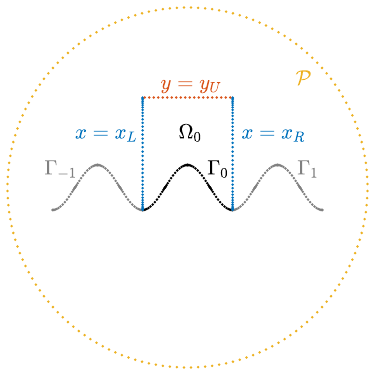}
    \caption{Unit cell with proxy surface.}
\end{subfigure}
\caption{\label{fig:unitcell} Illustration of (a) an example geometry and (b) the corresponding unit cell. (a) A portion of the domain $\Omega$ and boundary $\Gamma$ are shown for the cosine geometry. The point source $\vct{x}_0$ and the period $d$ are labeled. (b) The unit cell is shown with its lower boundary $\Gamma_0$ in black, left and right walls in blue, and upper wall in red. The unit cell domain $\Omega_0$ is labeled. The left and right neighbors, $\Gamma_{-1}$ and $\Gamma_1$, are shown in gray. The proxy circle $\mathcal{P}$ with charges $\{y_j\}_{j=1}^{N_{\mathrm{proxy}}}$ is shown in orange encompassing the unit cell and its two neighbors.}
\end{figure}

Within the unit cell, the solution is represented as the sum of an integral formulation using the two-dimensional, free-space Helmholtz Green's function $G_\omega(\vct{x},\vct{y}) = \frac{i}{4}H_0^{(1)}(\omega\|\vct{x}-\vct{y}\|)$
and an additional contribution to enforce quasiperiodicity. The additional contribution is a linear combination of basis functions defined as the potential due to proxy charges placed at $N_{\rm proxy}$ points on a circle $\mathcal{P}$ of radius
$R_{\rm proxy}\in\left[ \frac{3d}{2} , 2d\right]$ \cite{Cho:15}. An example of a proxy circle $\mathcal{P}$ with radius $R_{\rm proxy}=2$ is illustrated in Figure \ref{fig:unitcell}(b). Each point $\vct{y}_j$ on the proxy surface $\mathcal{P}$ has a corresponding proxy basis function $\phi_j(\vct{x})$ defined by
\begin{equation}
 \phi_j(\vct{x}) = \frac{\partial G_{\omega}}{\partial \vct{\nu}_j} (\vct{x},\vct{y}_j) + i\omega G_{\omega}(\vct{x},\vct{y}_j),
 \label{eq:proxy}
\end{equation}
where $\vct{\nu}_j$ is the outward-pointing normal vector at $\vct{y}_j$. 

Specifically, the solution in the unit cell is expressed as
\begin{equation}
\label{eq:lev1}
   u_{\kappa}(\vct{x}) = (\tilde{\mathcal{S}}_{\Gamma}\sigma)(\vct{x}) 
 + \sum_{j=1}^{N_{\rm proxy}}c_j\phi_j(\vct{x}),
\end{equation}
where 
\begin{equation} \label{eq:single} (\tilde{\mathcal{S}}_\Gamma\sigma)(\vct{x}) = \sum_{l=-1}^1 \alpha^{l}\int_{\Gamma_0}G_\omega(\vct{x},\vct{y}+l\vct{d})\sigma(\vct{y})\d s(\vct{y}).\\
\end{equation}
Note that the kernel of $(\tilde{\mathcal{S}}_\Gamma\sigma)$ includes contributions from the left ($\Gamma_{-1}$) and right ($\Gamma_1$) copies of the unit cell boundary, appropriately scaled by the Bloch phase $\alpha$, in addition to the contribution from the boundary $\Gamma_0$ included in the unit cell. The boundary charge density $\sigma$ and the coefficients for the proxy charges $\{c_j\}_{j=1}^{N_{\rm proxy}}$ are yet to be determined.

\begin{remark}
The proxy charges are the same as the equivalent charges used in the development of some kernel-independent FMMs \cite{GILLMAN2014poisson,YING2004kernel}. The idea is that the field due to anything outside of the proxy circle can be captured to the desired accuracy with appropriately weighted charges on the surface. The restriction of the integral operators in (\ref{eq:lev1}) to the unit cell and the left and right neighbors keeps the contributions that need to be scaled by the Bloch phase to a minimum. This is directly correlated to a reduced computational cost in constructing the linear system and evaluating the field.
\end{remark}

The representation (\ref{eq:lev1}) is valid inside the unit cell, but it does not satisfy the radiation condition above the unit cell. Thus, above the unit cell
(i.e. for points in the unit strip where $y>y_U$), the solution is given by a Rayleigh--Bloch expansion. This expansion is uniformly convergent \cite{bonnet1994guided} and satisfies the radiation condition (\ref{eq:ray}).
Specifically, for $\vct{x} = (x,y)$ in the unit strip where $y>y_U$, the solution is given by 
\begin{equation}u_{\kappa}(x,y) = \sum_{n\in\mathbb{Z}} a_n e^{{\rm i}\beta_nx}e^{{\rm i}k_n(y-y_U)}
 \label{eq:ral1}
\end{equation}
where $\beta_n$ and $k_n$ are the on-surface and vertical wavenumbers as defined in \eqref{eq:wavenumbers},
and the set $\{a_n\}$ are coefficients to be determined.

There are two big advantages to this representation over using the quasiperiodic Green's function as the kernel function. First, this representation is robust even at Wood's anomalies \cite{Cho:15}. Second, it can be split into terms dependent on and independent of the Bloch phase. This allows us to precompute large portions of the kernel, then scale them by the appropriate Bloch phase later \cite{Zhang21,zhang2022}.

\subsection{The linear system }
\label{sec:lin-sys}
The representations inside and outside the unit cell satisfy the partial differential equation in their respective regions, but they do not naturally enforce the boundary condition, quasiperiodicity of the solution, or continuity of the representations through the line $y = y_U$. 
Enforcing these conditions provides equations that can be solved to find the three sets of unknowns: $\sigma(x)$, $\{c_j\}_{j=1}^{N_{\rm proxy}}$, and $\{a_n\}$. 

The following boundary integral equation results from enforcing the 
Neumann boundary condition on $\Gamma_0$:
\begin{equation}
\label{eq:int}
-\frac{1}{2} \sigma(\vct{x})+(\tilde{\mathcal{D}}^*_{\Gamma}\sigma)(\vct{x}) 
 + \sum_{j=1}^{N_{\rm proxy}}c_j\partial_{\nu_{\vct{x}}}\phi_j(\vct{x}) = g(\vct{x}) \quad {\rm for } \  \vct{x}\in\Gamma_0,
\end{equation}
where 
\begin{equation}\label{eq:sumAdjoint}
(\tilde{\mathcal{D}}^*_{\Gamma}\sigma)(\vct{x}) =\sum_{l=-1}^1 \alpha^{l}\int_{\Gamma_0} \partial_{\nu_x}G_\omega(\vct{x},\vct{y}+l\vct{d})\sigma(\vct{y})\d s(\vct{y})
\end{equation}
is the adjoint double layer kernel and $g(\vct{x}) = -f_{\vct{\nu}}(\vct{x})$ from \eqref{eq:qp-rhs}. 

The equations that enforce quasiperiodicity of the solution and its flux are given by 
\begin{equation}
\label{eq:quasi}
\left.\left((\tilde{\mathcal{S}}_{\Gamma}\sigma)(\vct{x}) 
 + \sum_{j=1}^{N_{\rm proxy}}c_j\phi_j(\vct{x})\right)\right|_{x=x_L} = \alpha^{-1} \left.\left((\tilde{\mathcal{S}}_{\Gamma}\sigma)(\vct{x}) 
 + \sum_{j=1}^{N_{\rm proxy}}c_j\phi_j(\vct{x})\right)\right|_{x=x_R},
\end{equation}
and

\begin{equation}
\label{eq:quasiFlux}
\left.\left((\tilde{\mathcal{D}}^*_{\Gamma}\sigma)(\vct{x}) 
 + \sum_{j=1}^{N_{\rm proxy}}c_j\partial_{\nu_{\vct{x}}}\phi_j(\vct{x})\right)\right|_{x=x_L} = \alpha^{-1}\left.\left((\tilde{\mathcal{D}}^*_{\Gamma}\sigma)(\vct{x}) + \sum_{j=1}^{N_{\rm proxy}}c_j\partial_{\nu_{\vct{x}}}\phi_j(\vct{x})\right)\right|_{x=x_R},
\end{equation}
respectively. The last set of equations enforce continuity of the two representations (the integral representation (\ref{eq:lev1}) and the Rayleigh--Bloch expansion (\ref{eq:ral1})) and their derivatives through the top of the unit cell. The equations are given by
\begin{equation}
\label{eq:top1}
   \left.\left( (\tilde{\mathcal{S}}_{\Gamma}\sigma)(\vct{x}) 
 + \sum_{j=1}^{N_{\rm proxy}}c_j\phi_j(\vct{x})\right)\right|_{y=y_U} = 
\left.\left(\sum_{n\in\mathbb{Z}} a_n e^{{\rm i}\beta_nx}e^{{\rm i}k_n(y-y_U)} \right)\right|_{y=y_U}, 
\end{equation}
and
\begin{equation}
\label{eq:top2}
\left.\left(    (\tilde{\mathcal{D}}^*_{\Gamma}\sigma)(\vct{x}) 
 + \sum_{j=1}^{N_{\rm proxy}}c_j\partial_{\nu_{\vct{x}}}\phi_j(\vct{x})\right)\right|_{y=y_U} = 
\left.\left(\sum_{n\in\mathbb{Z}} a_n {\rm i}\beta_ne^{{\rm i}\beta_nx}e^{{\rm i}k_n(y-y_U)}\right)\right|_{y=y_U}. 
\end{equation}

Equations (\ref{eq:int}--\ref{eq:top2}) are continuous equations that cannot be solved directly; their solution needs to be approximated. We choose to approximate the integral operators via a Nystr\"om discretization. Let $\{\vct{x}_j\}_{j=1}^N$ denote the quadrature nodes on $\Gamma_0$ used in the Nystr\"om discretization. 
 In the same manner as \cite{Cho:15}, the quasiperiodicity of the solution and its flux are enforced at points that lie on Gaussian panels on the left and right walls of the unit cell. Let $M_w$ denote the number of points used to enforce quasiperiodicity. 
Lastly, the continuity of the integral representation and the Rayleigh--Bloch expansions is enforced at collection of  $M$ uniformly distributed points on the top of the unit cell. The Rayleigh--Bloch expansion is truncated at $n = \pm K$.

The rectangular linear system that arises from these choices has the form

\begin{equation}
 \mthree{\mtx{A}}{\mtx{B}}{\vct{0}}{\mtx{C}}{\mtx{Q}}{\vct{0}}{\mtx{Z}}{\mtx{V}}{\mtx{W}}\vthree{\vct{\sigma}}{\vct{c}}{\vct{a}}
 = \vthree{\vct{g}}{\vct{0}}{\vct{0}},
 \label{eq:bigblock}
\end{equation}
where  $\mtx{A}$ is a matrix of size $N\times N$,
$\mtx{B}$ is a matrix of size $N\times N_{\rm proxy}$, 
$\mtx{C}$ is a matrix of size $M_w\times N$,
$\mtx{Q}$ is a matrix of size $M_w\times N_{\rm proxy}$,
$\mtx{Z}$ is a matrix of size $(2M)\times N$,
$\mtx{V}$ is a matrix of size $(2M)\times N_{\rm proxy}$,
and $\mtx{W}$ is a matrix of size $(2M)\times (2K+1)$.
The three block row equations in the above system enforce the Neumann boundary condition, quasiperiodicity of the solution and its flux, and continuity of the integral representation and the Rayleigh--Bloch expansion, respectively.

Details on the specifics of the matrix entries are provided in \cite{Cho:15,Zhang21,zhang2022}.

\subsection{Representing the quasiperiodic Green's function}
\label{sec:gqp}
The boundary data $f_{\vct{\nu}}$ in the differential equation (\ref{eq:quasi-prob}) is the normal derivative of the quasiperiodic Green's function evaluated on the boundary $\Gamma$. The periodizing scheme described earlier in this section can be used to construct the quasiperiodic Green's function. The normal derivative of it can then be used to evaluate $f_{\vct{\nu}}$.

\begin{figure}
    \centering
    \includegraphics{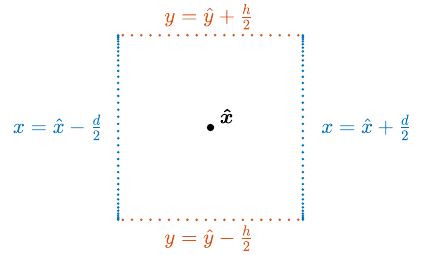}
    \caption{\label{fig:gqp-unitcell} Illustration of the unit cell used for constructing the quasiperiodic Green's function given a source point $\hat{\vct{x}}$. The left and right walls are shown in blue, top and bottom walls in red. The width of the unit cell is the period $d$ and the height is the parameter $h$.}
\end{figure}

Let $\hat{\vct{x}} = (\hat{x},\hat{y})$ denote the location of a source point. A different unit cell is defined as a rectangle of width $d$ and a chosen height $h$, as shown in Figure~\ref{fig:gqp-unitcell}. The unit cell is placed such that the point source is centered in the rectangle. The construction of this unit cell is like the one from section \ref{sec:geom} except there is no geometry to create the lower bound, so a bottom wall is created. The unit cell's walls are discretized as before, i.e. $M_w$ Gaussian nodes are placed on the left and right walls $\left(\text{at } x=\hat{x}-\frac{d}{2} \text{ and } x=\hat{x}+\frac{d}{2}\right)$, and $M$ uniformly distributed points are placed on the top and bottom walls $\left(\text{at } y=\hat{y}+\frac{h}{2} \text{ and } y=\hat{y}-\frac{h}{2}\right)$ each.

For $\vct{x}$ inside the unit cell, the quasiperiodic Green's function, i.e.\ the potential due to a line of point sources at locations $\hat{\vct{x}} + l\vct{d}$, $l \in \mathbb{Z}$, and Bloch phase $\alpha$, is given by 

\begin{equation}
\label{eq:gqp}
\psi(\vct{x},\hat{\vct{x}}) = \sum_{l=-1}^1\alpha^l G_\omega(\vct{x},\hat{\vct{x}}+l\vct{d}) +\sum_{j=1}^P\hat{c}_j \phi_j(\vct{x}).
\end{equation}

For $\vct{x}=(x,y)$ above and below the unit cell, the field is given by a Rayleigh--Bloch expansion \cite{bruno2014rapidly}. Thus for $y> \hat{y} + \frac{h}{2}$ and $y < \hat{y} - \frac{h}{2}$,
\begin{equation} \label{eq:gqp-rb}
    \psi(\vct{x},\hat{\vct{x}}) = \sum_{n\in\mathbb{Z}} \hat{a}_n e^{{\rm i}\beta_n x} e^{{\rm i}k_n\left|y-\left(\hat{y}+\frac{h}{2}\right)\right|},
\end{equation}
where $\beta_n$ and $k_n$ are as defined in \ref{eq:wavenumbers}. Let $\{\hat{a}_n\}$ denote the unknown coefficients of this expansion.

With these representations, it is possible to follow the example from the previous subsection to create a linear system. The big difference is that there is no integral operator and no boundary condition to satisfy. In other words, the system looks like (\ref{eq:bigblock}) except the first block row equation does not exist. The resulting linear system is given by
\begin{equation}
\label{eq:gqpsys}
    \begin{bmatrix}
        \mtx{Q} & \mtx{0}\\
        \mtx{V} & \mtx{W}
    \end{bmatrix}\begin{bmatrix}
        \hat{\vct{c}}\\
        \hat{\vct{a}}
    \end{bmatrix} = -\begin{bmatrix}
        \mtx{C}\\
        \mtx{Z}
    \end{bmatrix}
\end{equation}

where  
$\mtx{C}$ is a matrix of size $M_w\times 1$,
$\mtx{Q}$ is a matrix of size $M_w\times N_{\rm proxy}$,
$\mtx{Z}$ is a matrix of size $(4M)\times 1$,
$\mtx{V}$ is a matrix of size $(4M)\times N_{\rm proxy}$,
and $\mtx{W}$ is a matrix of size $(4M)\times (2K+1)$.

Here, the matrix entries of $\mtx{C}$
and $\mtx{Z}$ are the evaluations of 
$\sum_{l=-1}^1\alpha^l G_\omega(\vct{x},\hat{\vct{x}}+l\vct{d})$ and 
its flux on the left/right walls and the top/bottom of the unit cell, respectively. 
The other matrices are as defined in the last section.
In practice, the rectangular linear system in equation (\ref{eq:gqpsys}) is small and can quickly be solved using a backward stable application of the pseudoinverse.

\section{Efficient direct solution technique}
\label{sec:solver}

For each quadrature node $\kappa_j$ used for approximating the integral (\ref{eq:fb-integral}), the linear system (\ref{eq:bigblock}) must be solved. If the number
of points $N$ on the boundary is not large, it is possible to solve the system for each $\kappa_j$ via dense linear algebra. When $N$ is large (i.e. $N>10^4$ depending on computational resources), it is possible to build a fast direct solver like the ones presented in \cite{Zhang21,2013_Gillman,zhang2022}. This class of fast solvers is not necessary for the geometries of interest in this paper since they do not require a large number of discretization points to achieve the desired accuracy. This section presents the proposed solution technique which utilizes some ideas from  \cite{Zhang21,2013_Gillman,zhang2022} coupled with techniques that accelerate linear algebra for the matrix sizes under consideration.

Algorithm 2 provides an overview of the precomputation steps in the solution technique.

The proposed direct solver is based on the creation of a block solve for the linear system (\ref{eq:bigblock}). We chose to use a block matrix solve in the format
\begin{equation}
    \begin{array}{rl}
    \vct{b} & = \mtx{S}^{\dagger}\hat{\mtx{C}}\mtx{A}^{-1}\vct{g},\\
    \vct{\sigma} & = \mtx{A}^{-1}\vct{g} -\mtx{A}^{-1}\hat{\mtx{B}}\vct{b},
    \end{array}
    \label{eq:blocksolve}
\end{equation}
where $\vct{b} = \vtwo{\vct{c}}{\vct{a}}$, 
$\hat{\mtx{Q}} = \mtwo{\mtx{Q}}{\mtx{0}}{\mtx{V}}{\mtx{W}}$, 
$\hat{\mtx{C}} = \vtwo{\mtx{C}}{\mtx{Z}}$, $\hat{\mtx{B}} = \left[\mtx{B} \ \mtx{0}\right]$, and
\begin{equation}\label{eq:schur1}
\mtx{S} = \hat{\mtx{Q}}-\hat{\mtx{C}}\mtx{A}^{-1}\hat{\mtx{B}}
\end{equation}
denotes
the Schur complement matrix.
The matrix $\mtx{S}^{\dagger}$ denotes the pseudoinverse of $\mtx{S}$. 

To create a direct solver that is efficient for all $\kappa_j$, it 
is necessary to reuse as much computation. The first place for computational savings is in the evaluation of the linear system. 
Specifically, for all $\kappa_j$, the entries of all matrices except $\mtx{W}$ are made up of linear combinations of the same matrices, just scaled by appropriate amounts of the Bloch phase $\alpha$. Thus, all the unscaled matrices can be precomputed. 

The next place for acceleration is in the inversion of $\mtx{A}$ which must be done for each Bloch phase $\alpha$. The matrix $\mtx{A}$ can be written as
$$\mtx{A} = \mtx{A}_0 +\alpha\mtx{A}_{1} +\alpha^{-1}\mtx{A}_{-1},$$ 
where $\mtx{A}_0$ denotes the contribution due to the unit cell acting on itself, and $\mtx{A}_{\pm 1}$ denote the contributions from the left and right neighbor boundaries $\Gamma_{\pm 1}$ interacting with $\Gamma_0$, respectively. These correspond to $l=0,-1$, and $1$ in the sum of the boundary integral operators (\ref{eq:sumAdjoint}). 

The matrices $\mtx{A}_{-1}$ and $\mtx{A}_{1}$ are low rank because they predominantly correspond to interaction of well-separated points. In the low frequency regime, like the problems under consideration in this manuscript, the ranks are $O(1)$ \cite{1998_lowFrequency}. Let 
$\mtx{L}_{-1}\mtx{R}_{-1}$ and $\mtx{L}_{1}\mtx{R}_{1}$
 denote the low rank factorizations of $\mtx{A}_{-1}$ and $\mtx{A}_{1}$, respectively. Section \ref{sec:lowrank} presents the technique for creating the low rank factorizations. Let 
 $$\mtx{L} = \left[\mtx{L}_{-1} \ \mtx{L}_1\right] \quad {\rm and} \quad
\mtx{R} = \vtwo{\alpha^{-1}\mtx{R}_{-1}}{\alpha \mtx{R}_1}.$$ Then the inverse of $\mtx{A}$ can be approximated using these factorizations and the Woodbury formula,
\begin{equation}
\label{eq:woodbury}
\mtx{A}^{-1}\approx \left( \mtx{A}_0 +\mtx{L}\mtx{R} \right)^{-1} = \mtx{A}_0^{-1} - \mtx{A}_0^{-1}\mtx{L} (\mtx{I}+\mtx{R}\mtx{A}_0^{-1}\mtx{L})^{-1} \mtx{R}\mtx{A}_0^{-1}.
\end{equation} 
Note that the inverse $\mtx{A}_0^{-1}$ can be created once and used for all $\alpha$. For each $\alpha$, the sub-blocks of $\mtx{R}$ only need to be scaled appropriately. In practice, the matrix in (\ref{eq:woodbury}) is never created. Instead it is applied rapidly to any vector $\vct{f}$ using the formula 
\begin{equation}
\label{eq:woodburyapply}
 \left( \mtx{A}_0 +\mtx{L}\mtx{R} \right)^{-1}\vct{f} = \mtx{A}_0^{-1}\vct{f} - \mtx{A}_0^{-1}\mtx{L} (\mtx{I}+\mtx{R}\mtx{A}_0^{-1}\mtx{L})^{-1} \mtx{R}\mtx{A}_0^{-1}\vct{f},
\end{equation} 
where $\mtx{A}_0^{-1}\mtx{L}$ can be computed once and reused for all $\alpha$, and $\mtx{A}_0^{-1}\vct{f}$ is evaluated once per $\alpha$.

\subsection{Creating the low rank factorization of $\mtx{A}_{-1}$}
\label{sec:lowrank}
For simplicity of presentation, we present the technique for creating the low rank factorization of $\mtx{A}_{-1}$. We choose to create the low rank factorization using the interpolatory decomposition \cite{gu1996,lowrank} defined as follows.

\begin{definition}
 The \textit{interpolatory decomposition} of an $m\times n$ matrix $\mtx{M}$ that has rank $l$ is
the factorization 
$$ \mtx{M} \sim \mtx{P}\mtx{M}(J(1:l),:),$$
where $J$ is a vector of integers $j$ such that $1\leq j\leq m$, and $\mtx{P}$ is an $m\times l$ matrix that contains an $l \times l$ identity matrix.
Namely, $\mtx{P}(J(1:l),:) = \mtx{I}_l$. 
\end{definition}

The $l$ rows of the matrix $\mtx{M}$ that are picked via the interpolatory decomposition are called the \emph{skeleton rows}. When $\mtx{M}$ results from the discretization of an integral equation, each row corresponds to discretization points called the \emph{skeleton points}. For the problems under consideration in manuscript, a randomized linear algebra algorithm is used to find the rank $l$ to the user prescribed accuracy.

\begin{figure}
\centering
\begin{subfigure}[t]{0.37\textwidth}
    \centering
    \includegraphics{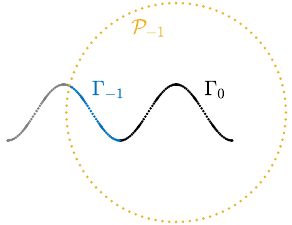}
    \caption{Full circle compression.}
\end{subfigure}
\hfill
\begin{subfigure}[t]{0.3\textwidth}
    \centering
    \includegraphics{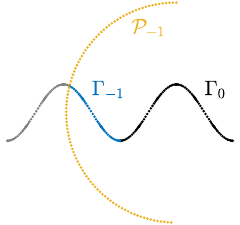}
    \caption{Directional compression.}
\end{subfigure}
\hfill
\begin{subfigure}[t]{.3\textwidth}
    \centering
    \includegraphics{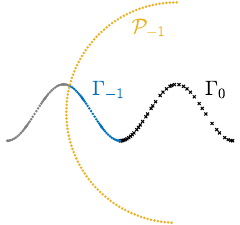}
    \caption{Skeleton points.}
\end{subfigure}
\caption{\label{fig:neigh}Illustrations of the proxy surfaces used to reduce the cost of creating the low rank factorization of $\mtx{A}_{-1}$. The unit cell boundary $\Gamma_0$ is shown in black, the near left neighbor boundary is shown in blue, and the far-left neighbor boundary is shown in gray. The proxy surface $\mathcal{P}_{-1}$ used in place of the far nodes is shown in orange. (a) The proxy surface is a full circle. (b) The proxy surface is a half-circle in the direction of the neighbor boundary. (c) The skeleton points found from the interpolatory decomposition index matrix $J(1:l)$ are denoted on the unit cell boundary with black crosses.}
\end{figure}

Creating the low rank factorization of $\mtx{A}_{-1}\in \mathbb{C}^{N\times N}$ has a computational cost that scales as $O(N^2l)$, where $l$ is the numerical rank of $\mtx{A}_{-1}$. Potential theory is useful to reduce this computational cost. Specifically, we place a collection of $n_{\rm proxy}$ proxy charges on a surface $\mathcal{P}_{-1}$ some distance ``far" from $\Gamma_0$. If the interpolatory decomposition can accurately capture the interaction of $\Gamma_0$ with $\mathcal{P}_{-1}$, then the interpolation matrix and skeleton can accurately capture the interaction of $\Gamma_0$ with any point left of $\mathcal{P}_{-1}$. Circles like the one illustrated in Figure \ref{fig:neigh}(a) are a common choice for the proxy surface \cite{Zhang21,2013_Gillman,zhang2022}. In the low frequency regime we are interested in, this choice can result in an artificially large rank. Instead, set $\mathcal{P}_{-1}$ to be a half circle as illustrated in Figure \ref{fig:neigh}(b).

Computing the interpolatory decomposition with $\mathcal{P}_{-1}$ alone is not sufficient as the interactions with 
$\Gamma_{-1}$ that lie inside $\mathcal{P}_{-1}$ (illustrated via the blue dots in Figure \ref{fig:neigh}) are not captured. 
The points on $\Gamma_{-1}$ that lie inside the proxy surface are called \textit{near points}. Let $I_{\rm near}$ denote the indices of the near points. 
Then the matrix $\mtx{A}_{\rm near} = \mtx{A}_{-1}(:,I_{\rm near})$ is the matrix that corresponds to the interaction of $\Gamma_0$ with the near points. 

Thus, to create $\mtx{L}_{-1}$ and $\mtx{R}_{-1}$, we apply the interpolatory factorization to 
$$\left[\mtx{A}_{\rm near} | \mtx{A}_{\rm proxy}\right],$$
where $\mtx{A}_{\rm proxy}$ is the matrix that captures interactions between $\Gamma_0$ and the proxy points on $\mathcal{P}_{-1}$. Let $J$ denote the index vector that results from the interpolatory decomposition. Then $\mtx{L}_{-1} = \mtx{P}$ resulting from the interpolatory decomposition and $\mtx{R}_{-1} = \mtx{A}_{-1}(J(1:l),:)$. Figure~\ref{fig:neigh}(c) 
illustrates the skeleton points obtained using the compression process using a half circle proxy surface.

\subsection{Corner compression}
\label{sec:corner}

\begin{figure}
\centering
\begin{subfigure}[t]{0.49\textwidth}
    \centering
    \includegraphics{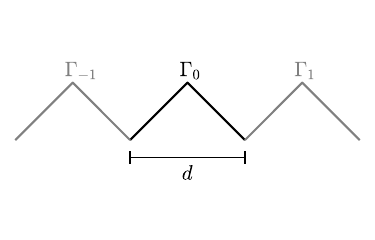}
    \caption{Staircase geometry.}
\end{subfigure}
\hfill
\begin{subfigure}[t]{0.49\textwidth}
    \centering
    \includegraphics{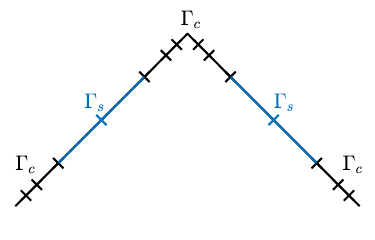}
    \caption{Corner refinement.}
\end{subfigure}
\caption{\label{fig:corner1} (a) Illustration of the partitioning of the staircase geometry with three periods. (b) Illustration of $\Gamma_0$ with the corner partitioning. $\Gamma_c$ in black denotes the part of the geometry near the corner that needs refinement, $\Gamma_s = \Gamma_0 \setminus \Gamma_c$ in blue. Shown are two levels of refinement into the corner. }
\end{figure}

When the boundary $\Gamma$ has corners, the density $\sigma$ has a singularity at each corner. The discretization technique must be modified to handle the singularity. An option is to use a quadrature that is specifically designed for corner singularities \cite{2016_corner,2019_corner2,2019_corner1}. More commonly, a panel-based quadrature is used, such as \cite{gen_quad} or \cite{helsing}. For simplicity, we consider a staircase geometry with three corners, shown in Figure \ref{fig:corner1}(a), to illustrate this method. Roughly speaking, the panel-based approach places a collection of panels on the boundary so that the part of the surface away from the corners are fully resolved. Then the panels near the corner are refined dyadically toward the corner to handle the singularities. Figure \ref{fig:corner1}(b) illustrates two levels of refinement in the corners of the staircase geometry. Although the latter refinement is easy to implement, many of the discretization points placed inside the corner are not necessary for accuracy purposes. This results in a linear system that is artificially large. One way to reduce the system size is to use the recursively compressed inverse preconditioning (RCIP) method \cite{helsing}. Since RCIP is difficult to extend to three-dimensional problems, we choose to use the corner compression method from \cite{GILLMAN2014}. 

The idea for this method is to partition the geometry according to proximity to the corner(s). Let $\Gamma_c$ denote the portion of the boundary in the corner(s) where the discretization is refined, and $\Gamma_s$ denote the remainder of the boundary. Figure \ref{fig:corner1}(b) illustrates the boundary partitioning for the staircase geometry that has one complete corner and two half corners at the left and right edges of the unit cell. Then the linear system resulting from discretizing a boundary integral equation on the geometry can be reordered as
\begin{equation}\mtwo{\mtx{A}_{cc}}{\mtx{A}_{cs}}{\mtx{A}_{sc}}{\mtx{A}_{ss}}\vtwo{\vct{q}_c}{\vct{q}_s} = \vtwo{\vct{f}_c}{\vct{f}_s}.
\label{eq:2by2}
\end{equation}

From potential theory, we know that the interactions between $\Gamma_c$ and $\Gamma_s $ are low rank. Let $l$ denote the rank of these interactions. Then the off-diagonal blocks in the system (\ref{eq:2by2}) can be factored as
\begin{equation}
\label{eq:lowrank}
\begin{array}{cccccccc}
\mtx{A}_{cs} &=& \mtx{U}_{c} & \mtx{B}_{cs} \\
n_{c} \times n_{s} && n_{c} \times l & l\times n_{s}
\end{array},
\quad\mbox{and}\quad
 \begin{array}{cccccccc}
\mtx{A}_{sc} &=&  \mtx{B}_{sc}& \mtx{V}^*_{c}  \\
n_{s} \times n_{c} && n_{s} \times l & l\times n_{c}
\end{array}.
\end{equation}

When (\ref{eq:lowrank}) holds, one can solve the following linear system 
instead of (\ref{eq:2by2}),
\begin{equation}
\label{eq:2x2_small}
\left[\begin{array}{cc}
\mtx{D}_{cc} & \mtx{B}_{cs} \\
\mtx{B}_{sc} & \mtx{A}_{ss}
\end{array}\right]\,
\left[\begin{array}{c}
\tilde{\vct{q}}_{c} \\ \vct{q}_{s}
\end{array}\right]
=
\left[\begin{array}{c}
\tilde{\vct{f}}_{c} \\ \vct{f}_{s}
\end{array}\right],
\end{equation}
where
\begin{equation}
\label{eq:formulas}
\mtx{D}_{cc} = \bigl(\mtx{V}_{c}^{*}\mtx{A}_{cc}^{-1}\mtx{U}_{c}\bigr)^{-1},
\qquad
\tilde{\vct{q}}_{c} = \mtx{V}_{c}^{*}\vct{q}_{c},
\qquad
\tilde{\vct{f}}_{c} = \mtx{D}_{cc}\mtx{V}_{c}^{*}\mtx{A}_{cc}^{-1}\vct{f}_{c}.
\end{equation}
The solution $\{\vct{q}_{c},\,\vct{q}_{s}\}$ of the larger system
(\ref{eq:2by2}) can be obtained from the solution $\{\tilde{\vct{q}}_{c},\,\vct{q}_{s}\}$
of the smaller system (\ref{eq:2x2_small}) via the formula
\begin{equation}\label{eq:full}
\vct{q}_{c} = \mtx{A}_{cc}^{-1}\,\mtx{U}_{c}\,\mtx{D}_{cc}\,\tilde{\vct{q}}_{c} + 
\bigl(\mtx{A}_{cc}^{-1} - \mtx{A}_{cc}^{-1}\,\mtx{U}_{c}\,\mtx{D}_{cc}\,\mtx{V}_{c}^{*}\,\mtx{A}_{cc}^{-1}\bigr)\,\vct{f}_{c}.
\end{equation}
The equivalence holds whenever the inverses in (\ref{eq:formulas}) exist. This method can be used to reduce the computational cost of inverting $\mtx{A}_0$ in the block system solver when there is refinement into a corner or multiple corners. 

In practice, we create the factorizations in (\ref{eq:lowrank}) via the interpolatory decomposition.  As as result, the matrices $\mtx{B}_{cs}$ and $\mtx{B}_{sc}$ are submatrices of $\mtx{A}_{cs}$ and $\mtx{A}_{sc}$, respectively, from equation (\ref{eq:2by2}).

\begin{remark}
When $\vct{f}_c$ is in the range of $\mtx{A}_{cs}$ it is possible to reduce \cite{GILLMAN2014} equation (\ref{eq:full}) to 
\begin{equation}\label{eq:range}\vct{q}_{c} = \mtx{A}_{cc}^{-1}\,\mtx{U}_{c}\,\mtx{D}_{cc}\,\tilde{\vct{q}}_{c}.\end{equation}
\end{remark}

Inverting $\mtx{A}_{cc}$ is not expensive because the matrix is small.

\subsubsection{A technique for compressing a single corner}
\label{sec:one-corner}	 	
\begin{figure}
\centering
\includegraphics{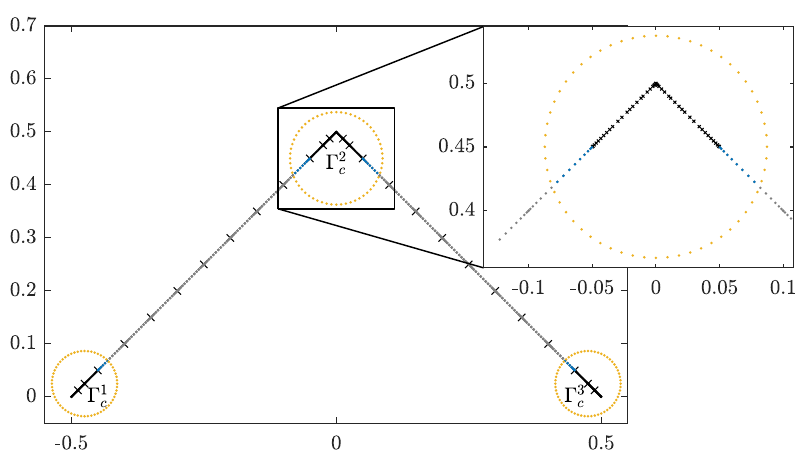}
\caption{\label{fig:corner}Illustration of a staircase geometry where the corner panels have two levels of refinement. There are three corners labeled $\Gamma_c^1$, $\Gamma_c^2$ and $\Gamma_c^3$. The refined panels at each corner are shown in black. The near part of the remaining boundary $\Gamma_s$ is shown in blue, and its far part in gray. The proxy surfaces used for compression of each corner are shown in orange. The inset shows the corner compression for $\Gamma_c^2$. The black crosses denote the skeleton points.} 
\end{figure}

The technique for creating the low rank factorizations associated with the corner compression follow the procedure outlined in section \ref{sec:lowrank}. First, proxy charges are placed on a circle centered around the corner boundary $\Gamma_c$. The near points correspond to points on $\Gamma_0$ that lie inside the circle but are on $\Gamma_s$. Figure \ref{fig:corner} illustrates the proxy circles and near points for the three corners of the staircase geometry $\Gamma_c^1$, $\Gamma_c^2$, and $\Gamma_c^3$. The interpolatory decomposition is used to create the low rank factorizations which \emph{skeletonize} the corner. The zoomed box inset in Figure \ref{fig:corner} illustrates the skeleton points for $\Gamma_c^2$ with black crosses.

\subsubsection{A technique for compressing multiple corners}
\label{sec:mult-corner}

When the boundary in the unit cell has multiple corners, we label the union of the corners $\Gamma_c$, i.e.
$\Gamma_c = \cup_{j=1}^{N_c} \Gamma_c^j$ where $N_c$ is the number of corners. For the $j^{\rm th}$ corner, the matrix $\mtx{A}_{c_js}$ defines the interaction of the corner with $\Gamma_s$. The low rank factorization of this matrix using the single-corner compression technique from section \ref{sec:one-corner} has factors $\mtx{P}_j$ and $\mtx{B}_j$. The matrix $\mtx{P}_j$ denotes the interpolation matrix and $\mtx{B}_j$ denotes the corresponding submatrix of $A_{c_js}$ containing the rows picked by the interpolatory decomposition.

The factorization defined for $\Gamma_c$ is given by putting these factorizations together to create the factors $\mtx{U}_c$ and $\mtx{B}_{cs}$ defined by
$$ \mtx{U}_c = \operatorname{diag} \left(\mtx{P}_1,\mtx{P}_2,\ldots,\mtx{P}_n\right), \qquad \mtx{B}_{cs} = \begin{bmatrix}
                    \mtx{B}_1 \\
                    \mtx{B}_2 \\
                    \vdots \\
                    \mtx{B}_n
                    \end{bmatrix}.
$$

The factorization of $\mtx{A}_{sc}$ is created via a similar procedure, taking into account that this matrix is in the transposed location.

\begin{algorithm}[ht]
\caption{Precomputation for the quasiperiodic solver}
\label{alg:precomp}
	\renewcommand{\algorithmicrequire}{\textbf{Input:}}
	\renewcommand{\algorithmicensure}{\textbf{Output:}}
\begin{algorithmic}[1]
    \Require wavenumber $\omega$ and 
    two-dimensional periodic boundary $\Gamma$.
    \Ensure precomputation of the Bloch phase independent portions of the solver.

    \State Compute all the unscaled matrices in $\mtx{B}$, $\mtx{C}$, $\mtx{Q}$, $\mtx{Z}$, and $\mtx{V}$.

    \State Compute $\mtx{A}_0$ and $\mtx{A}_0^{-1}$.

    \State Create the low rank factors $\mtx{A}_{-1} \approx \mtx{L}_{-1}\mtx{R}_{-1}$ and $\mtx{A}_1\approx \mtx{L}_{1}\mtx{R}_{1}$ via the technique in Section \ref{sec:lowrank}.

    \State With these factors, create the matrix $\mtx{L}=\left[\mtx{L}_{-1} \ \mtx{L}_1\right]$.

    \State Compute $\mtx{A}_{0}^{-1}\mtx{L}$.

    \Statex \textbf{For geometries with corners:}

    \State Skip the computation of $\mtx{A}_0$ and $\mtx{A}_0^{-1}$.

    \State Instead, create the interpolatory decomposition for the corners via the technique in Section \ref{sec:mult-corner} and build all operators needed for the linear system (\ref{eq:2x2_small}).

    \State Invert (\ref{eq:2x2_small}) directly.
    
    \State Compute $\mtx{A}_{cc}^{-1}$, $\mtx{A}_{cc}^{-1}\mtx{U}_c$, $\mtx{D}_{cc}$, and $\mtx{D}_{cc}\mtx{V}_c^*$.

    \State Use (\ref{eq:full}) to apply $\mtx{A}_0^{-1}$.
\end{algorithmic}
\end{algorithm}

\begin{algorithm}[ht]
	\renewcommand{\algorithmicrequire}{\textbf{Input:}}
	\renewcommand{\algorithmicensure}{\textbf{Output:}}
\caption{Application of the quasiperiodic solver}
\label{alg:solve}
\begin{algorithmic}[1]
    \Require the operators created via Algorithm 1, $\kappa$ which defines the Bloch phase $\alpha$, the location of point charge $\vct{x}_0$, and the boundary data $f_{\vct{\nu}}$ computed via Section \ref{sec:gqp}.
    \Ensure the approximation of the boundary charge density $\vct{\sigma}$, the proxy charge coefficients $\vct{c}$, and Rayleigh-Bloch expansion coefficients $\vct{a}$.
    \Statex The output allows for $u_\kappa$ to be evaluated at any point $\vct{x}=(x,y)$ via (\ref{eq:lev1}) if $y \leq y_U$ and via (\ref{eq:ral1}) if $y>y_U$.
    
    \State Set $\vct{g} = f_{\vct{\nu}}\big|_{\Gamma_0}$.

    \State Appropriately scale matrices with $\alpha$ to create $\mtx{R}$, $\mtx{B}$, $\mtx{C}$, $\mtx{Q}$, $\mtx{Z}$, and $\mtx{V}$.
    
    \State Compute $\mtx{W}$.
    
    \State Evaluate $\mtx{A}^{-1}\mtx{B}$ and $\mtx{A}^{-1}\vct{g}$ using the Woodbury formula (\ref{eq:woodburyapply}).
    
    \State Evaluate the Schur complement $\mtx{S}$ defined in (\ref{eq:schur1}).
    
    \State Create the SVD of $\mtx{S}$.
    
    \State Evaluate $\vct{b} = \mtx{S}^{\dagger}(\hat{\mtx{C}}\mtx{A}^{-1}\vct{g})$ using the SVD to apply $\mtx{S}^{\dagger}$.
    
    \State Evaluate $\vct{\sigma} = \mtx{A}^{-1}\vct{g} -\mtx{A}^{-1}\hat{\mtx{B}}\vct{b}$.

    \State \Return $\vct{\sigma}$, $\vct{c}$, and $\vct{a}$.
\end{algorithmic}
\end{algorithm}

\section{Numerical results}
\label{sec:numerics}

This section reports on the performance of the proposed solution technique. For two-dimensional problems, we expect the bulk of the computational gain to come from the choice of periodization method and the ability to reuse precomputed operators for all quasiperiodic solves required to approximate the Floquet--Bloch integral.

The numerical examples first illustrate the performance of the direct solution technique, allowing for the cost of the precomputation and the application of the solver to be easily separated. Then the solver is used to approximate the Floquet--Bloch integral (\ref{eq:fb-integral}), highlighting the performance of the full solution technique to find the scattered field due to a point source.

The accuracy of the direct solution technique is measured by the relative error of the quasiperiodic solution evaluated at a target point defined by
\begin{equation}\label{eq:rel-err}
    E_N= \left|\frac{u_N(\vct{x},\vct{x}_0) - u_{N_{\rm max}}(\vct{x},\vct{x}_0)}{u_{N_{\rm max}}(\vct{x},\vct{x}_0)}\right|,
\end{equation}
where $u_N$ is the approximate solution found with $N$ discretization points and $u_{N_{\rm max}}$ is the approximate solution with the largest number of discretization points. Descriptions of the distribution of discretization points are reported for each experiment. The relative error for the aperioidic solution due to a point source is reported following a similar format.

To illustrate the benefit of low rank linear algebra, we consider four different approaches to building the direct solver. For each approach, the $\alpha$-independent blocks are precomputed but the calculation of $\mtx{A}^{-1}$ differs. The approaches are listed and described below.
\begin{itemize}
    \item \textbf{Dense:} The matrix $\mtx{A}^{-1}$ is inverted directly in the solve. 
    \item \textbf{ID proxy:}  The precomputation includes the low rank factorizations of $\mtx{A}_{-1}$ and $\mtx{A}_1$ using a full proxy circle (see Figure \ref{fig:neigh}(a)), evaluating $\mtx{A}_0^{-1}$, and evaluating $\mtx{A}_0^{-1}\mtx{L}$. The Woodbury formula is used to apply $\mtx{A}^{-1}$ in the solve step. 
    \item \textbf{Half circle:} The precomputation includes the low rank factorizations of $\mtx{A}_{-1}$ and $\mtx{A}_1$ using a half circle as the proxy surface (see Figure \ref{fig:neigh}(b)), evaluating  $\mtx{A}_0^{-1}$ and evaluating $\mtx{A}_0^{-1}\mtx{L}$. The Woodbury formula is used to apply $\mtx{A}^{-1}$ in the solve step. 
    \item \textbf{Corner comp:}  The precomputation includes the low rank factorizations of $\mtx{A}_{-1}$ and $\mtx{A}_1$ using a half circle as the proxy surface (see Figure \ref{fig:neigh}(b)), the corner compression for constructing the action of $\mtx{A}_0$ and the matrix $\mtx{A}_0^{-1}\mtx{L}$. The Woodbury formula is used to apply $\mtx{A}^{-1}$ in the solve step.
\end{itemize}

The period of all geometries is $d = 1$. The boundaries are discretized via 16-node composite modified Gaussian quadrature \cite{KOLM_Rok,gen_quad} . For all experiments, $N$ denotes the total number of discretization points on the boundary in the unit cell, and the compression tolerance $\epsilon$ is set to $10^{-13}$.

The walls of the unit cell are set to be 1 unit in height. The number of Gauss--Legendre nodes for the left and right walls is set to $M_w = 240$, totaling 120 nodes on each wall. The number of equispaced nodes on the top wall is set to $M=60$. The truncation of the Rayleigh--Bloch expansion is set to $K=20$. The proxy circle $\mathcal{P}$ is set to have $N_{\rm proxy} = 160$ equispaced points. The radius of the proxy circle is $R_{\rm proxy} = 2d$.

For computing low rank factorizations, the number of points on the proxy surface is set to $n_{\rm proxy} = 100$. Let $R$ denote the radius of the smallest concentric circle containing the geometry segment of interest. Then the proxy surface used to create the low rank factorization is the concentric circle or half circle with a radius of $1.75R$.

The timings for using the quasiperioidic solver is broken into two stages: the precomputation stage and solve stage. The distribution of work for these stages are as defined in their corresponding algorithms: Algorithm \ref{alg:precomp} and Algorithm \ref{alg:solve}. When evaluating the Floquet--Bloch integral, the total time for creating the approximation is reported.

For all geometries and wavenumbers considered, the same contour deformation defined by  (\ref{eq:contour}) and illustrated in Figure~\ref{fig:fb-integral}(a) is used to approximate the Floquet--Bloch integral. We examine cases when the wavenumber $\omega$ is such that the contour passes closely between two coalescing singularities and when it does not. In the former case, quadrature nodes are clustered near the singularities according to (\ref{eq:expgrading}). When the contour is not near coalescing singularities, $N_{\kappa}$ periodic trapezoidal nodes are used to approximate the integral. All experiment specific parameters are stated in their corresponding sections.

All results were computed on a Dell Precision 7920 Tower, with two 12-core Intel Xeon Gold 6128 processors and 192 GB of memory. The algorithm is implemented in MATLAB. All low rank factorizations are created using the interpolatory decomposition routine in \texttt{librla}, a randomized linear algebra library implemented in MATLAB \cite{2026librla,librlaGIT}.

\subsection{Comparison with chunkIE}
\label{sec:chunkie}

In these first experiments, we compare the proposed solution technique to the solution technique available in an integral equation software package \texttt{chunkIE} for two geometries:  the cosine geometry shown in Figure \ref{fig:unitcell} and the staircase geometry shown in Figure \ref{fig:corner1}. \texttt{chunkIE} is a MATLAB package for solving two-dimensional partial differential equations using the boundary integral equation approach \cite{chunkiecode}. This package is one of the standard tools used by the boundary integral equations community. The quasiperiodic solution technique in \texttt{chunkIE} is based on using the quasiperioidic Green's function approach from \cite{agocs2024trapped} as implemented in \cite{agocs2025complexscaling}. In the reporting of the numerical results, we refer to this method as the \textbf{QP approach}. Comparisons of the two methods use the same boundary discretization points.

For the cosine geometry, the number of panels per period is $N_{\rm pan}$. For the staircase geometry, $0.5 N_{\rm pan}$ panels are placed on each segment of the staircase. The panel closest to each of the corners is refined dyadically with $N_{\rm ref}$ refinements. Since there are four panels in corners, the total number of panels on the staircase geometry is $N_{\rm pan}+4N_{\rm ref}$.

For all the experiments in this section, the source point is $\vct{x}_0 = (-0.2,0.35)$ and the point where the solution is evaluated at, i.e. the target point, is $\vct{x} = (0.3,0.25)$.

\subsubsection{Performance of the direct solver}

To illustrate the performance of the direct solution technique, we consider solving (\ref{eq:quasi-prob}) with $\omega=1.2$ and $\theta = \pi/5$. Recall that $\kappa = \omega\cos\theta$ and $\alpha = e^{i\kappa d}$.

\textit{Cosine geometry.} For the smooth cosine geometry, three approaches are used to build the direct solver: \textbf{Dense}, \textbf{ID proxy}, and \textbf{Half circle}. Table \ref{tab:cos-precomp-solve-times} reports the time in seconds for the precomputation and solve stages, and the relative error achieved using the different direct solver options with various $N_{\rm pan}$. The accuracy of the approximate solution is the same for all the direct solver options. While building the low rank factorizations makes the precomputation portion of the solver slower, the speed of applying the solver with the factorizations is much faster. In this example, with 100 or more panels, the factorization solver will be faster than the dense linear algebra solver any time the approximation of the integral (\ref{eq:fb-integral}) requires more than three solves. 

\renewcommand{\arraystretch}{1.2}
\begin{table}
    \centering
    \begin{tabular}{l l l l l l l l}
        \hline
        \multirow{2}{*}{\shortstack[l]{\textbf{Cosine}\\\textbf{Geometry}}} & $N_{\rm pan}$ & 4 & 8 & 20 & 40 & 100 & 200 \\
         & $N$ & 64 & 128 & 320 & 640 & 1600 & 3200 \\
        \hline
        Precomp & Dense & 0.095 & 0.16 & 0.38 & 0.75 & 2.4 & 5.9 \\
         & ID proxy & 0.11 & 0.21 & 0.47 & 0.87 & 2.5 & 7.0 \\
         & Half circle & 0.10 & 0.17 & 0.42 & 0.77 & 2.3  & 6.4 \\[2pt]
        Solve & Dense & 0.018 & 0.021 & 0.023 & 0.034 & 0.16 & 0.58 \\
         & ID proxy & 0.022 & 0.022 & 0.023  & 0.025  & 0.042 & 0.083 \\
         & Half circle & 0.020 & 0.019 & 0.022 & 0.026  & 0.040 & 0.069 \\[2pt]
        Error & Dense & 1.1e-10 & 1.0e-13 & 8.6e-14 & 8.1e-14 & 1.3e-13 & --- \\
         & ID proxy & 1.2e-10 & 9.9e-14 & 8.4e-14 & 5.0e-14 & 1.2e-13 & --- \\
         & Half circle & 1.1e-10 & 1.0e-13 & 1.7e-13 & 2.7e-13 & 2.3e-13 & --- \\
         \hline
    \end{tabular}
    \caption{Time in seconds taken by the precomputation and solve stages of the direct solver applied to the cosine geometry with $\omega = 1.2$ and $\theta = \pi/5$ for increasing number of panels $N_{\rm pan}$. The relative error, found using (\ref{eq:rel-err}), is also reported.}
    \label{tab:cos-precomp-solve-times}
\end{table}

\textit{Staircase geometry.}
For the experiments using the staircase geometry, 8 uniform panels are placed on each of the straight edges, i.e. $N_{\rm pan}=16$. This choice results in an approximation that is accurate to four digits. The panels lying closest to the corners are refined in a dyadic manner to gain further accuracy.
The number of refinements $N_{\rm ref}$ range from $0$ to $36$. In these experiments, we compare the \textbf{Dense}, \textbf{Half circle}, and \textbf{Corner comp} approaches. 

Table \ref{tab:stair-precomp-solve-times} reports the times taken by the precomputation and solve stages of the direct solvers. 
The precomputation time necessary for the Half circle and Corner comp methods is sometimes less than that required by the Dense approach because the full neighbor interaction matrices need not be constructed for the former two options. low rank linear algebra again reduces the solve time. As the staircase only has a few corners, compression of the self-interaction matrix has little impact on the solution time. For $N_{\rm ref}=12$ or greater, the Half circle approach only requires a single solve to be faster than the Dense approach.

\begin{table}
    \centering
    \begin{tabular}{l l l l l l l}
        \hline
        \multirow{2}{*}{\shortstack[l]{\textbf{Staircase}\\\textbf{Geometry}}} & $N_{\rm ref}$ & 0 & 6 & 12 & 24 & 36 \\
         & $N$ & 256 & 640 & 1024 & 1792 & 2560 \\
        \hline
        Precomp & Dense & 0.31 & 0.74 & 1.27 & 2.86 & 6.08 \\
         & Half circle & 0.35 & 0.82 & 1.25 & 2.82 & 5.32 \\
         & Corner comp & 0.40 & 0.96 & 1.68 & 3.04 & 5.87 \\[2pt]
        Solve & Dense & 0.022 & 0.033 & 0.059 & 0.19 & 0.36 \\
         & Half circle & 0.020 & 0.025 & 0.034 & 0.048 & 0.079 \\
         & Corner comp & 0.023 & 0.027 & 0.035 & 0.055 & 0.080 \\[2pt]
        Error & Dense & 2.9e-4 & 1.5e-5 & 9.8e-7 & 3.9e-9 & --- \\
         & Half circle & 2.7e-4 & 1.6e-5 & 1.0e-6 & 3.8e-9 & --- \\
         & Corner comp & 2.5e-4 & 1.5e-5 & 9.3e-7 & 3.9e-9 & --- \\
         \hline
    \end{tabular}
    \caption{Time in seconds taken by the precomputation and solve stages of the direct solver applied to the staircase geometry with $\omega = 1.2$ and $\theta = \pi/5$ for increasing number of refinements $N_{\rm ref}$ into each corner. The relative error, found using (\ref{eq:rel-err}), is also reported.}
    \label{tab:stair-precomp-solve-times}
\end{table}

Table \ref{tab:cornercompress} reports the size of the corner-compressed system, $N_{\rm compress}$, for different levels of refinements into each corner, $N_{\rm ref}$. As the level of refinement increases, the accuracy of the quasiperiodic solution improves (as shown in Table \ref{tab:stair-precomp-solve-times}), but the size of the compressed corner-to-corner interaction matrix $\mtx{D}_{cc}$ remains constant. This illustrates that the low rank factorization of the corner-to-edge interactions is doing a good job picking the ``important'' discretization points near the corners as skeleton points.

\begin{table}
\centering
\begin{tabular}{l l l}\hline
$N_{\rm ref}$ & $N$  & $N_{\rm compress}$\\ \hline
0 & 256 & 256 \\
6 & 640 & 292 \\
12 & 1014 & 292 \\
24 & 1702 & 292 \\
36 & 2560 & 292 \\
\hline
\end{tabular}
\caption{\label{tab:cornercompress}Table reporting the number of refinements into each corner $N_{\rm ref}$, the total number of
discretization points on the staircase $N$, and the size of compressed system $N_{\rm compress}$ when corner compression is used on the staircase geometry. The tolerance of the compression is held fixed at $\epsilon = 10^{-13}$.}
\end{table}

\subsubsection{Performance for the aperiodic boundary value problem}

Solving the aperiodic boundary value problem via the technique proposed in this manuscript requires approximating the Floquet--Bloch integral (\ref{eq:fb-integral}) using quadrature points and weights outlined in Section \ref{sec:flobloch}. This section 
illustrates the performance of the solution technique for evaluating that integral with three different choices of wavenumber: $\omega = 1.2$, $\omega = 2.4$, and $\omega = 0.01$. Any dashed entries in the timing results correspond to experiments that we could not conduct due to memory constraints.

\textit{Away from the branch points.}
When $\omega=1.2$, the contour deformation used in the Floquet--Bloch transform (\ref{eq:contour}) is away from the branch points. Thus, it is sufficient to use the periodic trapezoidal rule on the contour to achieve high accuracy. We set $N_\kappa=60$, meaning 60 quasiperiodic boundary value solves are required to evaluate $u$, the aperiodic scattered solution.

\begin{figure}
    \centering
    \begin{subfigure}[t]{\textwidth}
        \centering
        \includegraphics{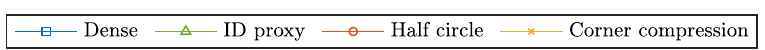}
        \vspace{1mm}
    \end{subfigure}
    \begin{subfigure}[t]{0.48\textwidth}
        \centering
        \includegraphics{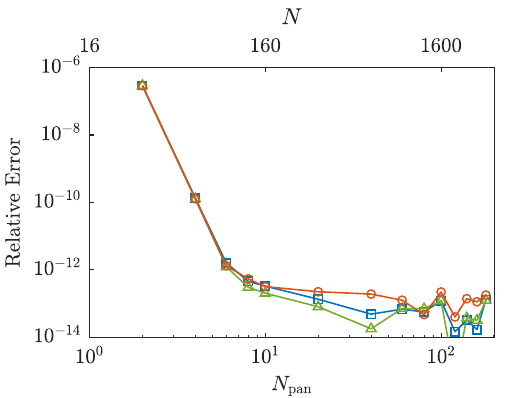}
        \caption{Cosine geometry.}
    \end{subfigure}
    \hfill
    \begin{subfigure}[t]{0.48\textwidth}
        \centering
        \includegraphics{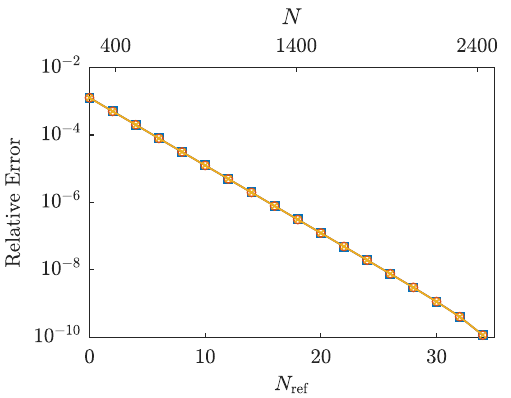}
        \caption{Staircase geometry.}
    \end{subfigure}
    \caption{The relative error in the Floquet--Bloch integral (i.e. the aperiodic solution) as a function of $N$, the total number of discretization points on the boundary, with $\omega = 1.2$ and $N_\kappa = 60$ for two geometries: (a) the cosine geometry, with the number of panels, $N_{\rm pan}$, on the boundary listed on the bottom axis, and (b) the staircase geometry, with the number of refinements, $N_{\rm ref}$, into each corner on the boundary listed on the bottom axis. }
    \label{fig:fullsoln-err}
\end{figure}

The first experiment involves the cosine geometry. Figure \ref{fig:fullsoln-err}(a) illustrates the relative error in approximating the solution due to a point source as the boundary is refined uniformly. The solution technique can achieve 13 digits of accuracy with roughly 10 panels (i.e.\ 160 discretization nodes) on $\Gamma_0$. Table \ref{tab:cos-fullsoln-times} reports the time in seconds required to evaluate the approximation. The benefit of using the low rank versions of the solver increases with the number of discretization points on the boundary. For $N_{\rm pan}=20$, our proposed method is about 3 times faster than the QP approach.

\begin{table}
    \centering
    \begin{tabular}{l l l l l l l}
        \multicolumn{7}{c}{\textbf{Cosine Geometry}} \\
        \hline
        $N_{\rm pan}$ & 4 & 8 & 20 & 40 & 100 & 200 \\
        $N$ & 64 & 128 & 320 & 640 & 1600 & 3200 \\
        \hline
        QP approach & 2.4 & 4.5 & 13 & 35 & 179 & --- \\
        Dense & 2.7 & 3.0 & 3.9 & 5.6 & 16 & 52 \\
        ID proxy & 2.8 & 3.3 & 4.4 & 6.2 & 13 & 23 \\
        Half circle & 2.7 & 3.1 & 4.0 & 5.8 & 11 & 20 \\
        \hline
    \end{tabular}
    \caption{The time in seconds for evaluating the Floquet--Bloch integral (\ref{eq:fb-integral}) for the cosine geometry with $\omega = 1.2$ and $N_\kappa=60$ Floquet--Bloch nodes using the QP approach, Dense, ID proxy and Half circle solvers for increasing number of panels $N_{\rm pan}$.}
    \label{tab:cos-fullsoln-times}
\end{table}

Next, we consider the staircase geometry. Figure \ref{fig:fullsoln-err}(b) illustrates the relative error for the solution due to a point source as the discretization is refined into the corner. Table 
\ref{tab:stair-fullsoln-times} reports the time in seconds required to evaluate the approximation. In this example,
the refinements are necessary to gain accuracy and using the low rank linear algebra is roughly twice as fast as using dense linear algebra. The proposed method is roughly 20 times faster than the QP approach for this problem when 9 digits of accuracy are desired.

\begin{table}
    \centering
    \begin{tabular}{l l l l l l}
        \multicolumn{6}{c}{\textbf{Staircase Geometry}} \\
        \hline
        $N_{\rm ref}$ & 0 & 6 & 12 & 24 & 36 \\
        $N$ & 256 & 640 & 1024 & 1792 & 2560 \\
        \hline
        QP approach & 11 & 38 & 80 & 235 & 539 \\
        Dense & 3.8 & 6.3 & 9.0 & 22 & 36 \\
        Half circle & 3.9 & 5.6 & 7.7 & 14 & 20 \\
        Corner comp & 4.2 & 6.4 & 8.5 & 13 & 19 \\
        \hline
    \end{tabular}
    \caption{The time in seconds for evaluating the Floquet--Bloch integral (\ref{eq:fb-integral}) for the staircase geometry with $\omega = 1.2$ and $N_\kappa=60$ Floquet--Bloch nodes using the QP approach, Dense, Half circle, and Corner comp solvers for increasing number of refinements, $N_{\rm ref}$ into each corner.}
    \label{tab:stair-fullsoln-times}
\end{table}

\textit{Near the branch points.}
When $\omega = 2.4$ and $\omega = 0.01$, the contour (\ref{eq:contour}) passes between two coalescing branch points and the refinement procedure defined in equation (\ref{eq:expgrading}) is necessary. For these experiments, the exponential grating parameter $b$ is set to $5$. For these experiments, the number of discretization points $N$ is held fixed and the number of quadrature nodes $N_{\kappa}$ for the Floquet--Bloch integral is varied.

For the cosine geometry, the number of panels is fixed to $N_{\rm pan} = 10$. Table \ref{tab:eg-times} shows the timing and error results of the full solution evaluation for various $N_\kappa$ nodes. Since many points are not required to resolve the geometry, the improvement from the proposed technique is small. This was expected from the results in the previous section. 

\begin{table}
    \centering
    \begin{tabular}{l l l l l l l l l l}
        \hline
          & $\boldsymbol{\omega}$ & \multicolumn{3}{l}{$\boldsymbol{2.4}$} & \multicolumn{3}{l}{$\boldsymbol{0.01}$} \\
         & $N_\kappa$ & 40 & 60 & 150 & 60 & 120 & 150 \\
        \hline
        \multicolumn{2}{l}{\textbf{Cosine geometry}} & & & & & & \\
        Time & QP approach & 5.77 & 6.26 & 15.0 & 4.89 & 8.41 & 10.3 \\
         & Half circle & 2.56 & 3.40 & 9.12 & 3.21 & 6.16 & 7.56 \\[2pt]
        Error & QP approach & 3.7e-5 & 9.6e-8 & --- & 2.7e-5 & 7.2e-9 & --- \\
         & Half circle & 3.6e-5 & 9.5e-8 & --- & 5.7e-5 & 9.6e-9 & --- \\
         \multicolumn{2}{l}{\textbf{Staircase geometry}} & & & & & & \\
        Time & QP approach & 122 & 504 & 1376 & 337 & 637 & 814 \\[2pt]
         & Corner comp & 13 & 15 & 32 & 11 & 19 & 23 \\
        Error & QP approach & 1.8e-7 & 7.5e-11 & --- & 2.8e-5 & 7.6e-9 & --- \\
         & Corner comp & 1.7e-7 & 7.0e-11 & --- & 6.0e-5 & 1.7e-8 & --- \\
         \hline
    \end{tabular}
    \caption{The time in seconds for evaluating the Floquet--Bloch integral (\ref{eq:fb-integral}) at wavenumbers $\omega = 2.4$ and $\omega = 0.01$ for the cosine and staircase geometries for increasing number of Floquet--Bloch nodes $N_\kappa$. The QP approach is compared to the Half circle and Corner comp solvers for cosine and staircase, respectively. The relative error, found using (\ref{eq:rel-err}), is also reported. $N_{\rm pan}=16$ for the cosine geometry, and $N_{\rm ref}=30$ for the staircase geometry.}
    \label{tab:eg-times}
\end{table}

For the staircase geometry, the number of panels and number of refinements are fixed to $N_{\rm pan}=16$ and $N_{\rm ref} = 30$. Table \ref{tab:eg-times} gives timing and error results of the full solution evaluation for various $N_\kappa$ nodes. For $\omega = 2.4$, the number of quadrature nodes needed on the contour (\ref{eq:expgrading}) to achieve
an 11-digit accurate approximation is $N_\kappa=60$. When $\omega = 0.01$, the contour is closer to the branch points. This means that additional refinement is needed to gain accuracy. In this experiment, $N_\kappa=120$ quadrature nodes were needed to 
obtain an approximation that is accurate to 8 digits. Table \ref{tab:eg-times} reports the timings for both the QP approach and the Corner comp approach. Since the case of $\omega=0.01$ requires many solves to achieve high accuracy, there is a $33$ times speedup using the Corner comp method instead of the QP approach.

\subsection{Complex geometries}
\label{sec:complex}
This section illustrates the performance of the proposed solver when applied to problems on two complex geometries. These geometries are non-symmetric and cornered boundaries. The geometries considered are a piecewise linear geometry with many corners, referred to as the dinosaur geometry, and a curved geometry with a single corner, referred to as the wave geometry.

For all experiments, the point source is located at $\vct{x}_0 = (0.4,0.45)$ and the target point is $\vct{x} = (-0.35,0.4)$. Three approaches are compared in this section: \textbf{Dense}, \textbf{Half circle}, and \textbf{Corner comp}.

\begin{figure}
    \centering
    \begin{subfigure}[t]{.49\textwidth}
        \centering
        \includegraphics{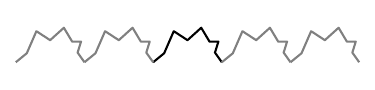}
        \caption{Dinosaur geometry.}
    \end{subfigure}
    \hfill
    \begin{subfigure}[t]{.49\textwidth}
        \centering
        \includegraphics{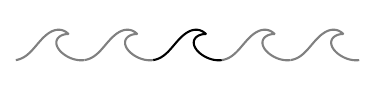}
        \caption{Wave geometry.}
    \end{subfigure}
    \caption{Illustration of the complex geometries. Five periods are illustrated. One period is highlighted in black.}
    \label{fig:complex-geoms}
\end{figure}

\textit{Dinosaur geometry.} The dinosaur geometry illustrated in Figure \ref{fig:complex-geoms}(a) is defined by connecting the vertices
$$  (-.50,.00),\ 
    (-.34,.13),\ 
    (-.20,.45),\ 
    ( .00,.32),\ 
    ( .20,.50),\ 
    ( .32,.30),\ 
    ( .45,.30),\ 
    ( .40,.14),\ 
    ( .50,.00)$$
with straight lines to form a piecewise linear boundary with nine corners per period, one at each vertex.

For the dinosaur geometry, take $N_{\rm pan}$ to be the number of panels per edge before refinement. 8 uniform panels are placed on each edge; i.e. $N_{\rm pan}=8$. Each corner is then dyadically refined $N_{\rm ref}$ times. Thus the total number of panels on one period of the geometry is $8N_{\rm pan}+16N_{\rm ref}$.

\textit{Wave geometry.} The wave geometry illustrated in Figure \ref{fig:complex-geoms}(b) is defined by a piecewise smooth curve with a single corner per period at the peak of the wave. The curve is defined peak-to-peak by a sixth-order B\'{e}zier curve with control points
$$
    \vct{p}_0 = (0.8,0.4),\  
    \vct{p}_1 = (0.5, 0.7),\ 
    \vct{p}_2 = (0.0,0.0),\ 
    \vct{p}_3 = (0.6,-0.3),\ 
    \vct{p}_4 = (-0.7,0.0),\ 
    \vct{p}_5 = (-0.5,0.4),\ 
    \vct{p}_6 = (-0.2,0.4).
$$
The parametric curve is given by
$$\vct{b}(t) = \sum_{k=0}^6 \binom{6}{k} \left(1-t\right)^{\left(6-k\right)} t^{k} \vct{p}_k, \quad 0\leq t \leq 1.$$

For the wave geometry, take $N_{\rm pan}$ to be the number of panels on the unit cell boundary before refinement. 8 panels are placed uniformly on one period of the geometry; i.e $N_{\rm pan}=8$. The panels near the corner are then dyadically refined $N_{\rm ref}$ times. Thus the total number of panels on one period of the geometry is $N_{\rm pan}+2N_{\rm ref}$.

\subsubsection{Performance of the direct solver}
In this section, we approximate the solution to (\ref{eq:quasi-prob}) with $\omega=1.2$ and $\theta = \pi/5$.

\textit{Dinosaur geometry.}
Table \ref{tab:dinowave-precomp-solve-times} reports the times for the precomputation and solve stages of the direct solver for the three approaches using the dinosaur geometry. The number of corners in this geometry is large enough that compressing corners for the self-interaction matrix decreases the time for applying the solver. For $N_{\rm ref}=24$ levels of refinement, after two solves the Corner comp method is faster than the Dense solver. After seven solves, the Corner comp method is faster than the Half circle method.

\textit{Wave geometry.}
Table \ref{tab:dinowave-precomp-solve-times} also reports the times for the precomputation and solve stages of the direct solver using the three approaches for the wave geometry. Compression of the self-interaction matrix sees less benefit and can even slow the solver because there is only one corner. There is modest speed up for using the Half circle method. For $N_{\rm ref}=24$ levels of refinement, the Half circle approach requires at least four solves to be faster than the Dense solver.

\begin{table}
    \centering
    \begin{tabular}{l l l l l l l l l l l l}
        \hline
        & & \multicolumn{5}{l}{\textbf{Dinosaur Geometry}} & \multicolumn{5}{l}{\textbf{Wave Geometry}} \\
         & $N_{\rm ref}$ & 0 & 6 & 12 & 24 & 36 & 0 & 6 & 12 & 24 & 30 \\
         & $N$ & 1024 & 2560 & 4096 & 7168 & 10240 & 128 & 320 & 512 & 896 & 1088 \\
        \hline
        Precomp & Dense & 1.35 & 5.01 & 9.99 & 24.05 & 43.94 & 0.20 & 0.45 & 0.68 & 1.22 & 1.51 \\
         & Half circle & 1.42 & 5.05 & 10.64 & 31.29 & 61.22 & 0.22 & 0.47 & 0.76 & 1.29 & 1.61\\
         & Corner comp & 1.81 & 5.61 & 11.24 & 32.07 & 61.19 & 0.25 & 0.53 & 0.80 & 1.41 & 1.75 \\[2pt]
        Solve & Dense & 0.057 & 0.42 & 1.11 & 4.06 & 9.51 & 0.020 & 0.024 & 0.029 & 0.047 & 0.061 \\
         & Half circle & 0.031 & 0.068 & 0.13 & 0.32 & 0.57 & 0.022 & 0.022 & 0.025 & 0.028 & 0.030 \\
         & Corner comp & 0.049 & 0.076 & 0.11 & 0.20 & 0.33 & 0.012 & 0.023 & 0.026 & 0.033 & 0.035 \\[2pt]
        Error & Dense & 4.4e-4 & 3.2e-5 & 2.1e-6 & 9.9e-9 & --- & 1.4e-3 & 1.3e-4 & 1.2e-5 & 8.9e-8 & --- \\
         & Half circle & 4.3e-4 & 3.1e-5 & 2.1e-6 & 9.4e-9 & --- & 1.4e-3 & 1.3e-4 & 1.2e-5 & 8.9e-8 & --- \\
         & Corner comp & 4.1e-4 & 3.0e-5 & 2.1e-6 & 9.9e-9 & --- & 1.4e-3 & 1.3e-4 & 1.2e-5 & 8.9e-8 & --- \\
         \hline
    \end{tabular}
    \caption{Time in seconds taken by the precomputation and solve stages of the direct solver applied to the dinosaur and wave geometries with $\omega = 1.2$ and $\theta = \pi/5$ for increasing number of refinements $N_{\rm ref}$ into each corner. The relative error, found using (\ref{eq:rel-err}), is also reported.}
    \label{tab:dinowave-precomp-solve-times}
\end{table}

\subsubsection{Performance for the aperiodic boundary value problem}
This section reports on the performance of the proposed technique for approximating the Floquet--Bloch integral (\ref{eq:fb-integral}) with the dinosaur and wave geometries for three choices of wavenumber: $\omega = 1.2$. $\omega = 2.4$, and $\omega = 0.01$.

\textit{Away from the branch points.}
For $\omega = 1.2$, $N_\kappa = 60$ periodic trapezoidal nodes are used to approximate the Floquet--Bloch integral. Three methods are compared here: \textbf{Dense}, \textbf{Half circle}, and \textbf{Corner comp}.

\begin{figure}
    \centering
    \begin{subfigure}[t]{\textwidth}
        \centering
        \includegraphics{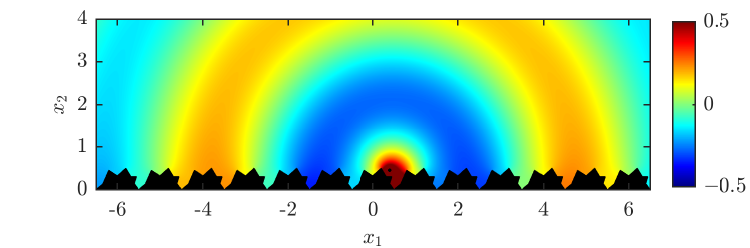}
        \caption{Dinosaur geometry.}
    \end{subfigure}
    \begin{subfigure}[t]{\textwidth}
        \centering
        \includegraphics{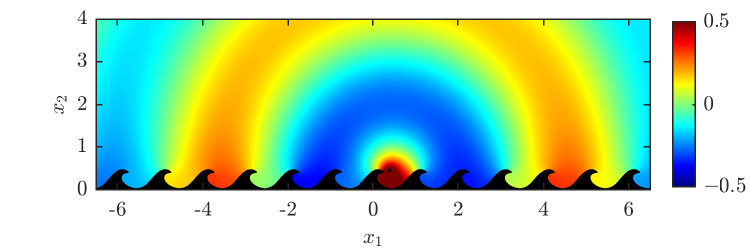}
        \caption{Wave geometry.}
    \end{subfigure}
    \caption{Real part of the total field, $u_{\rm tot} = u + u_i$, due to a point source $\vct{x}_0 = (0.4,0.45)$ with wavenumber $\omega = 1.2$, run with $N_{\rm ref}=24$. Top: dinosaur geometry accurate to 8 digits. Bottom: wave geometry accurate to 7 digits.}
    \label{fig:fullsoln-targs-dino-wave}
\end{figure}

Table \ref{tab:dinowave-fullsoln-times} reports the timing and error results for the dinosaur geometry. This geometry has many corners, so there is a benefit to using the corner compression technique. The Corner comp method is 4 times faster than the Dense approach when 8 digits of accuracy are desired. The real part of the total field due to a point source is shown in Figure \ref{fig:fullsoln-targs-dino-wave}(a) for many periods of the dinosaur geometry.

Table \ref{tab:dinowave-fullsoln-times} also reports the timing and error results for the wave geometry. This boundary does not require many discretization points to resolve, so there is not a significant benefit to using the compression techniques. With $N_{\rm ref} = 24$ levels of refinement, 7 digits of accuracy are achieved, confirming that curved geometries with corners are amenable to the proposed methods. The total field due to a point source is shown in Figure \ref{fig:fullsoln-targs-dino-wave}(b) for many periods of the wave geometry.

\begin{table}
    \centering
    \begin{tabular}{l l l l l l l l l l l l}
        \hline
         & & \multicolumn{5}{l}{\textbf{Dinosaur Geometry}} & \multicolumn{5}{l}{\textbf{Wave Geometry}} \\
        \multirow{2}{*}{\shortstack[l]{\textbf{}\\\textbf{}}} & $N_{\rm ref}$ & 0 & 6 & 12 & 24 & 36 & 0 & 6 & 12 & 24 & 30 \\
         & $N$ & 1024 & 2560 & 4096 & 7168 & 10240 & 128 & 320 & 512 & 896 & 1088 \\
        \hline
        Time & Dense & 10 & 38 & 89 & 282 & 625 & 3.5 & 4.9 & 6.3 & 9.4 & 11 \\
         & Half circle & 8.6 & 19 & 34 & 78 & 133 & 3.4 & 4.7 & 5.9 & 8.1 & 9.1 \\
         & Corner comp & 10 & 20 & 34 & 72 & 115 & 3.6 & 4.7 & 5.7 & 8.0 & 9.4 \\[2pt]
        Error & Dense & 1.1e-3 & 6.7e-5 & 4.4e-6 & 1.9e-8 & --- & 2.5e-3 & 2.2e-4 & 2.0e-5 & 1.5e-7 & --- \\
         & Half circle & 1.1e-3 & 6.7e-5 & 4.4e-6 & 1.9e-8 & --- & 2.5e-3 & 2.2e-4 & 2.0e-5 & 1.5e-7 & --- \\
         & Corner comp & 1.1e-3 & 6.7e-5 & 4.3e-6 & 1.9e-8 & --- & 2.5e-3 & 2.2e-4 & 2.0e-5 & 1.5e-7 & --- \\
         \hline
    \end{tabular}
    \caption{Runtime time in seconds for evaluating the Floquet--Bloch integral (\ref{eq:fb-integral}) for the dinosaur and wave geometries with $\omega = 1.2$ and $N_\kappa=60$ Floquet--Bloch nodes using the Dense, Half circle, and Corner comp solvers, as a function of the number of refinements into each corner $N_{\rm ref}$. The relative error, found using (\ref{eq:rel-err}), is also reported.}
    \label{tab:dinowave-fullsoln-times}
\end{table}

\textit{Near the branch points.}
For $\omega = 2.4$ and $\omega = 0.01$, clustered quadrature nodes are required to approximate the Floquet--Bloch integral accurately. The refinement procedure in equation (\ref{eq:expgrading}) is again used with $b=5$. The number of discretization points $N$ is held fix for the dinosaur and wave geometries and the number fo quadrature nodes $N_{\kappa}$ are varied.

For the dinosaur geometry, the number of panels is set to $N_{\rm pan} = 8$ and the number of refinements to $N_{\rm ref} = 30$. Table \ref{tab:eg-times-dinowave} shows the timing and error results for the aperiodic solution technique for various $N_{\kappa}$ nodes. 10 digits of accuracy are achieved with $N_{\kappa} = 60$ for $\omega = 2.4$. Additional nodes are needed to maintain accuracy for $\omega = 0.01$. $N_{\kappa} = 120$ gives 8 digits of accuracy. When 150 solves are used, the Corner comp method sees a 9 times speedup over the Dense solver.

For the wave geometry, the number of panels is set to $N_{\rm pan} = 8$ and the number of refinements to $N_{\rm ref} = 30$. Table \ref{tab:eg-times-dinowave} shows the timing and error results for the aperiodic solution technique for various $N_{\kappa}$ nodes. In line with the previous results, the wave geometry does not require many discretization points and thus does not benefit much from the Corner comp method. Similar accuracy is achieved as it was for the dinosaur geometry. This again shows that our methods work as expected on a curved, cornered geometry.

\begin{table}
    \centering
    \begin{tabular}{l l l l l l l l l l}
        \hline
          & $\boldsymbol{\omega}$ & \multicolumn{3}{l}{$\boldsymbol{2.4}$} & \multicolumn{3}{l}{$\boldsymbol{0.01}$} \\
         & $N_\kappa$ & 40 & 60 & 150 & 60 & 120 & 150 \\
        \hline
        \multicolumn{2}{l}{\textbf{Dinosaur geometry}} & & & & & & \\
        Time & Dense & 264 & 382 & 854 & 340 & 654 & 802 \\
         & Corner comp & 70 & 80 & 127 & 62 & 75 & 85 \\[2pt]
        Error & Dense & 2.7e-7 & 5.2e-10 & --- & 5.9e-5 & 7.7e-8 & --- \\
         & Corner comp & 2.7e-7 & 5.1e-10 & --- & 5.9e-5 & 7.7e-8 & --- \\
         \multicolumn{2}{l}{\textbf{Wave geometry}} & & & & & & \\
        Time & Dense & 6.6 & 9.0 & 20 & 6.6 & 13 & 16 \\[2pt]
         & Corner comp & 6.1 & 7.8 & 17 & 6.3 & 10 & 13 \\
        Error & Dense & 3.2e-7 & 4.2e-10 & --- & 5.9e-5 & 6.7e-8 & --- \\
         & Corner comp & 3.2e-7 & 4.2e-10 & --- & 5.9e-5 & 6.7e-8 & --- \\
         \hline
    \end{tabular}
    \caption{The time in seconds for evaluating the Floquet--Bloch integral (\ref{eq:fb-integral}) at wavenumbers $\omega = 2.4$ and $\omega = 0.01$ for the dinosaur and wave geometries for increasing number of Floquet--Bloch nodes $N_\kappa$. The Dense solver is compared to the Corner comp solver. The relative error, found using (\ref{eq:rel-err}), is also reported. The number of refinements is set to $N_{\rm ref}=30$.}
    \label{tab:eg-times-dinowave}
\end{table}

\section{Concluding remarks}
\label{sec:conclusion}

This manuscript presents an efficient direct solution technique for quasiperiodic scattering problems. The solver is built for the periodizing scheme originally presented in \cite{Cho:15}. Unlike previous fast solver papers \cite{Zhang21,zhang2022,2013_Gillman}, the solver presented in this paper excels for geometries that do not a large number of discretization points to get high accuracy. Instead, it uses ideas from the construction of fast direct solvers to compress interactions and create an efficient solver. The benefit of the solution technique is observed when approximating the Floquet--Bloch integral which requires many solves. For the problems in this paper, the proposed solution technique observes 20-30 times speed up for evaluating the Floquet--Bloch integral for a geometry with corners when compared with a solution technique based on using the quasiperiodic Green's function as the kernel function. 

This manuscript is the first to directly compare using the quasiperiodic Green's function with the periodization technique presented in section \ref{sec:boundInt}. We illustrated that simply changing to this representation and exploiting the ability to reuse the matrices results in massive speedups. This manuscript also introduced a directional compression technique that keeps the cost of creating the low rank factorizations with neighbors low while keeping the corresponding ranks close to optimal. This manuscript is also the first presentation of corner compression used to accelerate a quasiperiodic scattering solver. 

The authors are looking into extending this method to three-dimensional problems.
For three-dimensional, bi-periodic structures, the inverse Floquet--Bloch transform becomes a nested integral, requiring an even larger number of quasiperiodic solves to evaluate it accurately. Each solve involves a three-dimensional periodic surface, which increases the system size. Putting the bulk of the computational cost into a precomputation step will therefore become crucial for developing a computationally viable method. The extension of the direct solver to three-dimensional problems is not trivial. In particular, the directional compression techniques used in this paper do not directly extend to three-dimensional problems. Additionally, the evaluation of the quasiperiodic Green's function to create the boundary data $f_{\vct{\nu}}$ is more expensive because the periodizing points live on a sphere, a much larger number of points are needed on this surface than for the two-dimensional problem, and a matrix of size roughly the number of points on the sphere needs to be inverted for each new Bloch phase in the Floquet--Bloch integral.

\section{Acknowledgments}
The work by A. Gillman is supported by the National Science Foundation (DMS-2110886), and a 
Knut and Alice Wallenberg Foundation Grant. A. Gillman conducted a portion of this work while visiting the Institut Mittag-Leffler.

\bibliographystyle{abbrv} 
\bibliography{refs2}






\end{document}